\documentclass[11pt]{article}
\usepackage{mathrsfs}
\usepackage{upgreek}
\usepackage{float}
\usepackage{amssymb,amsmath,amsthm,amsfonts}
\usepackage[dvips]{graphicx}
\usepackage{epsfig, rotate}
\usepackage{pdfpages}

\usepackage[margin=1.2in]{geometry}

\numberwithin{equation}{section}

\newtheorem{theorem}{Theorem}[section]

\newtheorem{lemma}{Lemma}[section]

\usepackage{comment}

\def\eop{\hfill $\Box$ \medskip}

\newcommand\td{\overset{d}{\longrightarrow}}

\begin{document}
 \thispagestyle{empty}

\null\vspace{30pt}

\noindent {\bf \Large Limiting Distributions of the Likelihood Ratio
Test Statistics \\
for Independence of Normal Random Vectors}

\vspace{20pt}

\noindent {\bf Mingyue Hu$^{a, b}$,~~~ Yongcheng Qi$^{a}$}



\vspace{20pt}

\noindent $^a$Department of Mathematics and Statistics, University
of Minnesota Duluth, 1117 University Drive, Duluth, MN 55812, USA

\noindent $^b$School of Mathematics and Statistics, Southwest
University, 400715 Chongqing, China


\vspace{10pt}





\vspace{20pt}

\noindent\textbf{Abstract.} Consider the likelihood ratio test (LRT)
statistics for the independence of sub-vectors from a $p$-variate
normal random vector. We are devoted to deriving the limiting
distributions of the LRT statistics based on a random sample of size
$n$. It is well known that the limit is chi-square distribution when
the dimension of the data or the number of the parameters are fixed.
In a recent work by Qi, Wang and Zhang (Ann Inst Stat Math (2019)
71: 911--946), it was shown that the LRT statistics are
asymptotically normal under condition that the lengths of the normal
random sub-vectors are relatively balanced if the dimension $p$ goes
to infinity with the sample size $n$.  In this paper, we investigate
the limiting distributions of the LRT statistic under general
conditions. We find out all types of limiting distributions and
obtain the necessary and sufficient conditions for the LRT statistic
to converge to a normal distribution when $p$ goes to infinity. We
also investigate the limiting distribution of the adjusted LRT test
statistic proposed in Qi, Wang and Zhang (2019). Moreover, we
present simulation results to compare the performance of  classical
chi-square approximation, normal and non-normal approximation to the
LRT statistics, chi-square approximation to the adjusted test
statistic, and some other test statistics.

\vspace{10pt}




\vspace{10pt}


\newpage

\section{Introduction}\label{intro}

In classical theory for inference, many statistical methods have
been developed to test parametric hypotheses in the last few
decades. One of the most popular methods is the so-called likelihood
ratio test (LRT). It is known that the distribution functions of LRT
statistics can be well approximated by a chi-square distribution
under certain regularity conditions when the dimension of the data
or the number of the parameters of interest are fixed.  This means
that one does not have to estimate the variance of the test
statistics based on the likelihood ratio.

Many modern data sets such as financial data and modern
manufacturing data are  high-dimensional. Classical methods may not
be adequate for high dimensional data anymore, especially when the
dimension of data is relatively large compared with the sample size.
Some recent papers have investigated the limiting distribution of
the LRT statistics concerning the dependence structures of the
multivariate normal distributions.  It turns out that the chi-square
approximation fails while the dimension of the data increases with
the sample size. Instead, the normal approximation to the LRT
statistics works well under high dimension setting.  See, e.g., Bai
et al.~\cite{BJYZ}, Jiang et al.~\cite{JJY12}, Jiang and
Yang~\cite{JY13}, Jiang and Qi~\cite{JQ15},  Qi et al.~\cite{QWZ19},
Dette and D\"ornemann~\cite{DD2020}, and Guo and Qi~\cite{GuoQi21}.
Several approaches other than likelihood ratio method have been
developed in the literature; See, e.g., Schott~\cite{Schott2001,
Schott2005, Schott2007}, Ledoit and Wolf~\cite{LW02}, Bao et
al.~\cite{Bao2017}, Chen et al.~\cite{Chen2010},  Srivastava and
Reid~\cite{SR12}, Jiang et al~\cite{JBZ13},  Li et
al.~\cite{Li2017}, and Bodnar et al.~\cite{BDP2019}. A very recent
work by D\"ornemann~\cite{Dornemann2022} also established the
central limit theorems for some LRT statistics under non-normality.

In this paper, we consider a $p$-variate normal random vector and
study the limiting distributions of the LRT for testing the
independence of its grouped components based on a random sample of
size $n$. For a $p$-dimensional multivariate normal distribution
with mean vector $\mu$ and covariance matrix $\Sigma$, denoted by
$N_p(\mu, \Sigma)$, we partition its $p$ components into $k$ subsets
and test whether the $k$ sub-vectors are mutually independent, or
equivalently, we test whether the covariance matrix $\Sigma$ is
block diagonal.  In this paper, both $p$ and $k$ can depend on $n$
and diverge with the sample size. On the condition that the lengths
of the $k$ sub-vectors are relatively balanced, Qi et
al.~\cite{QWZ19} proved the asymptotic normality and proposed an
adjusted test statistic that usually has a chi-square limit when the
dimension $p$ goes to infinity with the sample size.

The aim of this paper is to give a complete description for the
limiting distributions of the LRT statistic for independence for
multivariate normal random vectors. We obtain all possible limiting
distributions and give the necessary and sufficient conditions for
the central limit theorem.  We also investigate the limiting
distributions of the adjusted test statistic proposed by Qi et
al.~\cite{QWZ19}, which performs better than the normal
approximation and chi-square approximation to the LRT statistics in
general.

The rest of the paper is organized as follows. In
Section~\ref{main}, we present our main results and establish the
necessary and sufficient conditions for the central limit theorem as
well as the conditions for non-normal limits. In Section~\ref{sim},
we present some simulation results to compare the performance of
four methods including the chi-square approximation to the LRT
statistic, the normal and non-normal approximation to the LRT
statistic and the chi-square approximation to the adjusted LRT
statistic proposed in Qi et al.~\cite{QWZ19}. In
Section~\ref{lemma}, we present some preliminary lemmas which are
used in Section~\ref{proof} to prove the main results of the paper.


\section{Main results}\label{main}

Let $\chi^2_f$ denote the random variables following chi-square
distribution with $f$ degrees of freedom and $N(0,1)$ the standard
normal variables.

For $k\geq 2$, let $q_1,\cdots, q_k$ be $k$ positive integers.
Denote $p=q_1+q_2+...+q_k$ and let
$$\Sigma=(\Sigma_{ij})_{1\leq i,j\le k}$$
be a positive definite matrix, where $\Sigma_{ij}$ is a $q_i\times
q_j$ sub-matrix. Assume $\xi_i$ is a $q_i$-dimensional normal random
vector for each $1\leq i \leq k$, and the $p$-dimensional random
vector $(\xi'_1, \cdots, \xi'_k)'$ has a multivariate normal
distribution $N_p(\mu, \Sigma)$. We are interested in testing the
independence of $k$ random vectors $\xi_1, \cdots, \xi_k$, or
equivalently the following hypotheses
 \begin{equation}\label{t1}
H_0: \Sigma_{ij}=0~\text{for all}~~1\leq i<j \leq k ~~~\text{vs}~~~
H_1: H_0~\text{is not true.}
 \end{equation}

Assume that $\mathbf{x}_1,\cdots, \mathbf{x}_n$ are $n$ independent
and identically distributed random vectors from distribution
$N_p(\mu, \Sigma)$. Define
\[
\mathbf{A}=\sum_{i=1}^{n}(\mathbf{x}_i-\mathbf{\bar
x})(\mathbf{x}_i-\mathbf{\bar x})',~~~ \mathbf{\bar
x}=\frac{1}{n}\sum_{i=1}^{n}\mathbf{x}_i,
\] and partition $\mathbf
A$ as follows
\[
\left(
  \begin{array}{cccc}
    \mathbf{A}_{11} & \mathbf{A}_{12} & \cdots & \mathbf{A}_{1k}\\
    \mathbf{A}_{21} & \mathbf{A}_{22} & \cdots & \mathbf{A}_{2k}\\
    \vdots & \cdots & \cdots & \vdots\\
    \mathbf{A}_{k1} & \mathbf{A}_{k2} & \cdots & \mathbf{A}_{kk}\\
  \end{array}
\right),
\]
where $\mathbf{A}_{ij}$ is a $q_i\times q_j$ matrix. According to
Theorem 11.2.3 from Muirhead~\cite{muirhead82},
the likelihood ratio statistic for testing \eqref{t1} is given by\\
\begin{equation}\label{equ1}
\Lambda_n=\frac{|\mathbf{A}|^\frac{n}{2}}{\prod_{i=1}^k|\mathbf{A}_{ii}|^\frac{n}{2}}:=(W_n)^{\frac{n}{2}}.
\end{equation}

Note that the likelihood ratio statistic $\Lambda_n$ is well defined
only if $p<n$.  When $p\ge n$, the determinant $|\mathbf{A}|$ is
zero since $\textbf{A}$ is singular.  Therefore, we can only
consider the case $p<n$ in the paper.

We introduce some notations before we give the main results. Let $g$
be any function defined over $(0,\infty)$. For integers $q$ and $n$
with $1\leq q < n$, define
\begin{equation}\label{Delta}
\Delta_{g,\,n,\,q}(t)=\sum_{i=1}^{q}g(\frac{n-i}{2}+t), ~~
t>-\frac{n-q}{2}.
\end{equation}
Let $\Gamma(x)$ denote the Gamma function, given by
\[
\Gamma(x)=\int^{\infty}_{0}t^{x-1}e^{-t}dt, ~x>0,
\]
and define the digamma function $\psi$
\begin{equation}\label{digamma}
\psi(x)=\frac{d\log
\Gamma(x)}{dx}=\frac{\Gamma'(x)}{\Gamma(x)},~x>0.
\end{equation}

\begin{theorem}\label{thm1}
Assume $p=p_n$ satisfy $2\leq p<n$ and $p_n\rightarrow\infty$ as
$n\rightarrow\infty$. Assume that $q_1,q_2,\cdots,q_k$ are $k$
positive integers such that $p=\sum_{i=1}^{k}q_i$, where $k=k_n$ may
depend on $n$. Set $q_{\max}=\max\{q_1,\cdots,q_k\}$ and assume
$n-q_{\max}\rightarrow \infty$ as $n\rightarrow\infty$. $\Lambda_n$
is the Wilks likelihood ratio statistic defined in $\eqref{equ1}$.
Then, under the null hypothesis in $\eqref{t1}$
\begin{equation}\label{T0}
T_0:=\frac{-2\log \Lambda_n-\mu _n}{\sigma _n}\td N(0,1)
\end{equation}
as $n\rightarrow \infty$, where
\begin{equation}\label{Mean}
\mu_n=-n\big(\Delta_{\psi,\,n,\,p}(0)-\sum_{j=1}^{k}\Delta_{\psi,\,n,\,q_j}(0)\big),
\end{equation}
\begin{equation}\label{Var}
\sigma_n^2=2n^2\big(\sum_{j=1}^{p}\frac{1}{n-j}-\sum_{i=1}^{k}
\sum_{j=1}^{q_j}\frac{1}{n-j}\big)+2n^2\big(b(n,p)-\sum_{i=1}^kb(n,
q_i)\big),
\end{equation}
\begin{equation}\label{bnq}
b(n,q)=\sum^q_{j=1}\frac{1}{(n-j)^2}, ~~1\le q<n,
\end{equation}
and symbol $\td$ denotes convergence in distribution.
\end{theorem}

Now we consider the situation when $n-q_{\max}$ is bounded. In this
case,  both $n-p$ and $p-q_{\max}$ are bounded because
$n-p+p-q_{\max}=n-q_{\max}$.  The following theorem gives non-normal
limits for $-2\log\Lambda_n$ when both $n-p$ and $p-q_{\max}$ are
fixed integers.

\begin{theorem}\label{thm2}
Assume $p=p_n$ satisfy $2\leq p<n$. Assume $q_1,q_2,\cdots,q_k$ are
$k$ positive integers such that $p=\sum_{i=1}^{k}q_i$, where $k=k_n$
is an integer that may depend on $n$. Set
$q_{\max}=\max\{q_1,\cdots,q_k\}$ and assume $p-q_{\max}=r$ and
$n-p=v$ for some fixed integers $r\ge 1$ and $v\ge 1$ for all large
$n$. Then, under the null hypothesis in \eqref{t1} we have
\begin{equation}\label{non-normal}
\frac{-2\log\Lambda_n-n\log n}{n}\td -\sum^{r+v-1}_{j=v} \log Y_j,
\end{equation}
where $Y_j$, $j\ge 1$ are independent random variables and the $Y_j$
has a chi-square distribution with $j$ degrees of freedom.
\end{theorem}

\noindent\textbf{Remark 1}. The classical likelihood method
considers the case when both $p$ and $k$ are fixed integers.  Assume
that $q_1, q_2, \cdots, q_k$ are fixed for all large $n$, then
\begin{equation}\label{classic}
-2\rho_n \log \Lambda_n\overset{d}\longrightarrow \chi_f^2,
\end{equation}
where
\begin{equation}\label{fn}
f  =  \frac{1}{2} \Big( p^2 - \sum_{i=1}^k q_i^2 \Big),
\end{equation}
\begin{equation}\label{rho}
\rho_n  =  1 - \frac{2 \Big( p^3 - \sum_{i=1}^k q_i^3 \Big) + 9
\Big( p^2 - \sum_{i=1}^k q_i^2 \Big)} {6n \Big( p^2 - \sum_{i=1}^k
q_i^2 \Big)};
\end{equation}
See, e.g., Theorem 11.2.5 in Muirhead~\cite{muirhead82}.

\noindent\textbf{Remark 2}.  Under conditions  that $q_{\max}\le
\delta p$ for some $\delta\in (0,1)$ and $p\to\infty$ as
$n\to\infty$, Qi et al.~\cite{QWZ19} established the following
central limit theorem for $-2\log \Lambda_n$
\[
T_1:=\frac{-2\log \Lambda_n -
\bar\mu_n}{\tau_n}\overset{d}\longrightarrow N(0, 1)
\]
as $n\to\infty$, where
\begin{equation}\label{mu1}
\bar\mu_n=n\sum_{i=1}^k(q_i-n+\frac32)\log(1-\frac{q_i}{n})-n(p-n+\frac32)\log(1-\frac{p}{n})+\frac{n}{3}\big(b(n,p)-\sum_{i=1}^kb(n,
q_i)\big),
\end{equation}
\[
\tau_n^2=2n^2\big(\sum_{i=1}^k\log(1-\frac{q_i}n)-\log(1-\frac{p}n)\big)+2n^2\big(b(n,p)-\sum_{i=1}^kb(n,
q_i)\big).
\]

\noindent\textbf{Remark 3.}  Theorem~\ref{thm1} is still true if
$\mu_n$ is replaced by $\bar\mu_n$ defined in \eqref{mu1} and
$\sigma_n^2$ is replaced by $\bar\sigma_n^2$
\begin{equation}\label{bar-Var}
\bar\sigma_n^2=2n^2\big(\sum_{j=1}^{p}\frac{1}{n-j}-\sum_{i=1}^{k}
\sum_{j=1}^{q_j}\frac{1}{n-j}\big).
\end{equation}
In fact, we have
\begin{equation}\label{mean-var}
\lim_{n\to\infty}\frac{\mu_n-\bar\mu_n}{\sigma_n}=0~~\mbox{
and}~~\lim_{n\to\infty}\frac{\bar\sigma_n^2}{\sigma_n^2}=1
\end{equation}
if $n-q_{\max}\to\infty$ as $n\to\infty$.  The proof of
\eqref{mean-var} is given in Lemma~\ref{Variance-Appr}.  In
practice, $\mu_n$ and $\sigma_n^2$ should be used since they give
better approximation to the mean and variance of $-2\log\Lambda_n$,
and this selection can achieve a better accuracy for the normal
approximation even when $n-q_{\max}$ is not very large.

\noindent\textbf{Remark 4.} The distribution of the random variable
on the right-hand side of \eqref{non-normal} is non-normal. This can
be verified by using the moment-generating functions.  The moment
generating function of $\sum^{r+v-1}_{j=v} \log Y_j$ is equal to
\[
2^{rt}\prod^{r+v-1}_{j=v}\frac{\Gamma(j/2+t)}{\Gamma(j/2)},
\]
which cannot be equal to $\exp(\mu t+ \sigma^2t^2/2)$, the moment
generating function of a normal random variable with a mean $\mu$
and variance $\sigma^2$. Otherwise, after taking the logarithm, it
implies that the third derivative of
$\sum^{r+v-1}_{j=v}\log\frac{\Gamma(j/2+t)}{\Gamma(j/2)}$ is
identically equal to zero.  We can show this cannot be true by using
some properties of the gamma function.  The details are omitted
here.

Now we are ready to establish the necessary and sufficient
conditions under which the central limit theorem holds for
$-2\log\Lambda_n$.

\begin{theorem}\label{thm3} There exist constants $a_n\in R$ and
$b_n>0$ such that under the null hypothesis in \eqref{t1}
\begin{equation}\label{clt}
\frac{-2\log\Lambda_n-a_n}{b_n}\td N(0,1)
\end{equation}
if and only if $p_n\to\infty$ and $n-q_{\max}\to\infty$ as
$n\to\infty$.
\end{theorem}

\vspace{10pt}

From \eqref{classic}, \eqref{T0} and \eqref{non-normal}, we have
only three different types of limiting distributions for
$-2\log\Lambda_n$, including chi-square distributions, normal
distributions and distributions of  linear combinations of
logarithmic chi-square random variables.  Theoretically,  for any
$p$ and $q_i$'s one can use one of the three limiting distributions
to approach the distribution of $-2\log\Lambda_n$ when $n$ is large.
Since the convergence in \eqref{classic}, \eqref{T0} or
\eqref{non-normal} does not provide clear cutoff values for $p$ and
$n-q_{\max}$, it may be difficult to select a limiting distribution
in practice even if $n$ is very large.

Qi et al.~\cite{QWZ19} proposed an adjusted log-likelihood ratio
test statistic (ALRT) which can be approximated by a chi-square
distribution when \eqref{classic} or \eqref{T0} holds. The ALRT is a
linear function of $-2\log\Lambda_n$ defined as
\begin{equation}\label{zn}
Z_n=(-2\log\Lambda_n)\sqrt{\frac{2f_n}{\sigma_n^2}}+f_n-\mu_n\sqrt{\frac{2f_n}{\sigma_n^2}}
\end{equation}
with $f_n$, $\mu_n$ and $\sigma_n^2$ being defined in \eqref{fn},
\eqref{Mean} and \eqref{Var}, respectively.  We have the following
result on chi-square approximation to the distribution of $Z_n$.

\begin{theorem}\label{thm4}
Let $p=p_n$ be a sequence of integers with $2\le p<n$. Assume
$k=k_n$ is also a sequence of positive integers, and $q_1,\cdots,
q_k$ are $k$ positive integers such that
 $p=\sum_{i=1}^k q_i$. Assume $n-q_{\max}\to\infty$ as $n\to\infty$. Then, under the null hypothesis in \eqref{t1}, we have
\begin{equation}\label{chisquare}
\lim_{n\to\infty}\sup_{-\infty<x<\infty}|P(Z_n\le
x)-P(\chi^2_{f_n}\le x)|=0.
\end{equation}
\end{theorem}

We note that the chi-square approximation in \eqref{chisquare} does
not impose any restriction on dimension $p$ and the number of
degrees of freedom of the chi-square distribution changes with $n$.

\section{Simulation study}\label{sim}

\subsection{Comparison of LRT related tests}\label{sub1}

In this subsection, we carry out a finite-sample simulation study to
compare the performance of three different approaches to the
likelihood ratio test statistic $-2\log\Lambda_n$ under condition
$n-q_{\max}\to\infty$ under the null hypothesis in \eqref{t1},
including the classical chi-square approximation \eqref{classic},
the normal approximation \eqref{T0}, and the adjusted chi-square
approximation in \eqref{chisquare}.  The three approaches are
denoted  by ``Chisq", ``CLT"  and ``ALRT", respectively. We will
also compare the performance of the normal approximation \eqref{T0},
the adjusted chi-square approximation in \eqref{chisquare} and the
non-normal approximation given in \eqref{non-normal} (denoted by
``LogChi"), and the comparison is made under conditions that both
$r=p-q_{\max}$ and $v=n-p$ are fixed integers.


For all four approaches, we will demonstrate how well the proposed
limiting distributions fit the histograms of the four test
statistics. From $\eqref{MGF}$ in Lemma~\ref{lemma1}, the
moment-generating function of $\log W_n$ is distribution-free under
the null hypothesis in $\eqref{t1}$. Since the four test statistics
are functions of $\Lambda_n$ and hence they are also functions of
$\log W_n$ from $\eqref{equ1}$, they are distribution-free as well
under the null hypothesis in $\eqref{t1}$. Therefore, the underlying
distribution in our study is set to be a multivariate normal
distribution with independent standard normal components.

In our simulation study, we choose sample size $n=101$.  For each of
selected combinations of $p$ and $q_i$'s under the regimes of the
chi-square and normal limiting distributions, we repeat the sampling
for $10,000$ times and obtain 10,000 replicates for the four test
statistics given in \eqref{classic}, \eqref{T0}, \eqref{zn}, and
\eqref{non-normal}. We plot the histogram for each test statistic
and its corresponding theoretical density function in one graph.

In the simulation study we consider the following four cases.

\begin{figure}[p!]
\centering
\includegraphics[height=.2\textheight, width=0.9\textwidth]{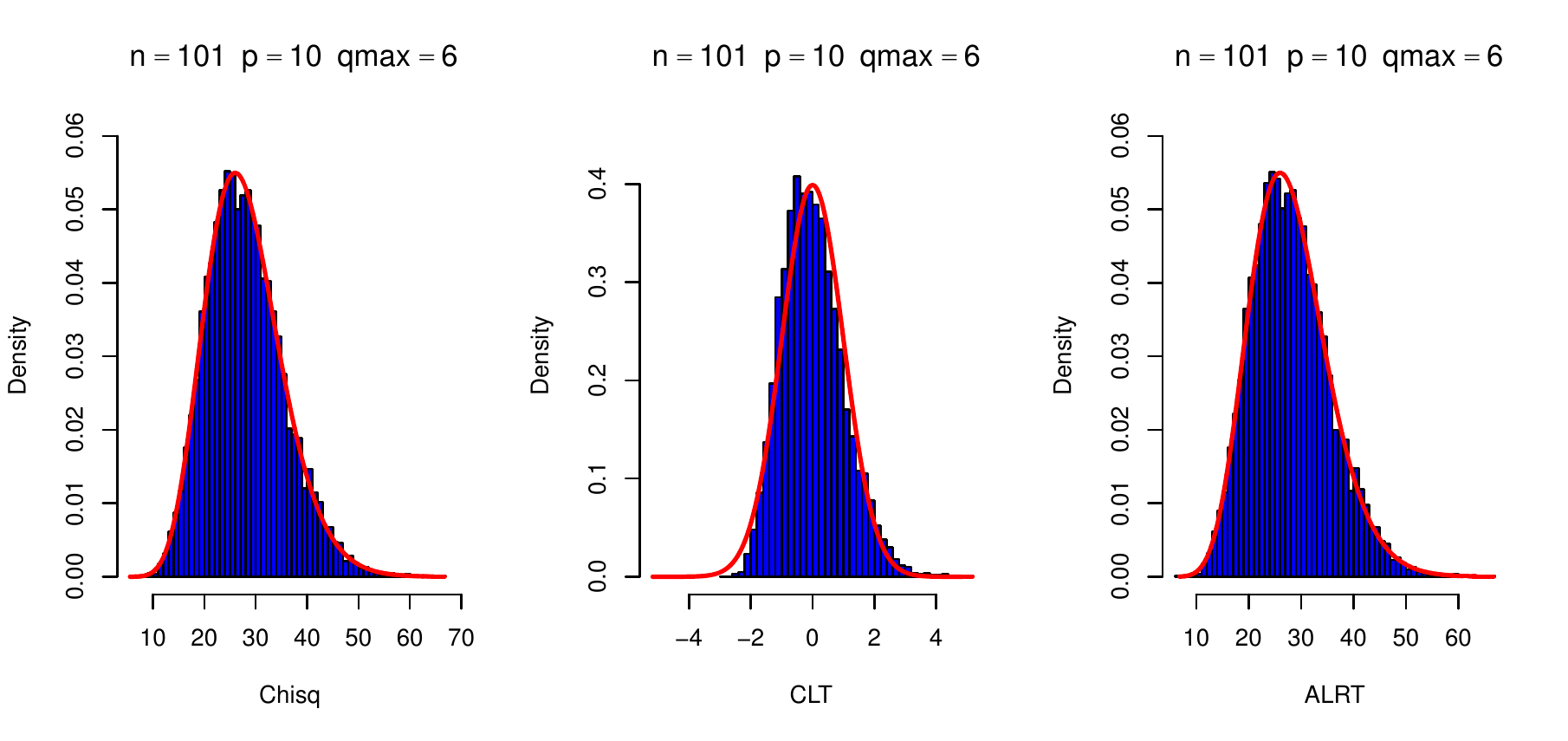}

\includegraphics[height=.2\textheight, width=0.9\textwidth]{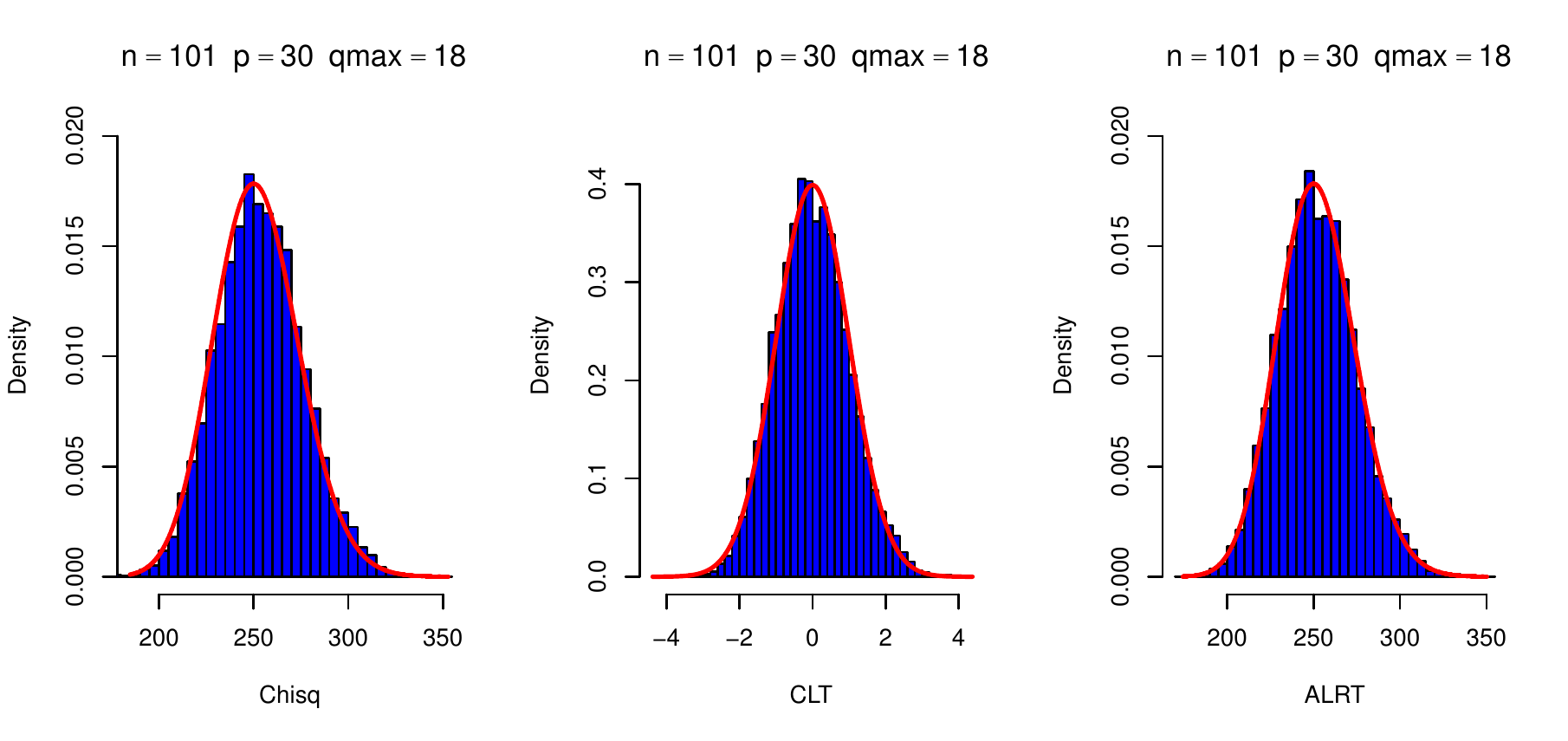}

\includegraphics[height=.2\textheight, width=0.9\textwidth]{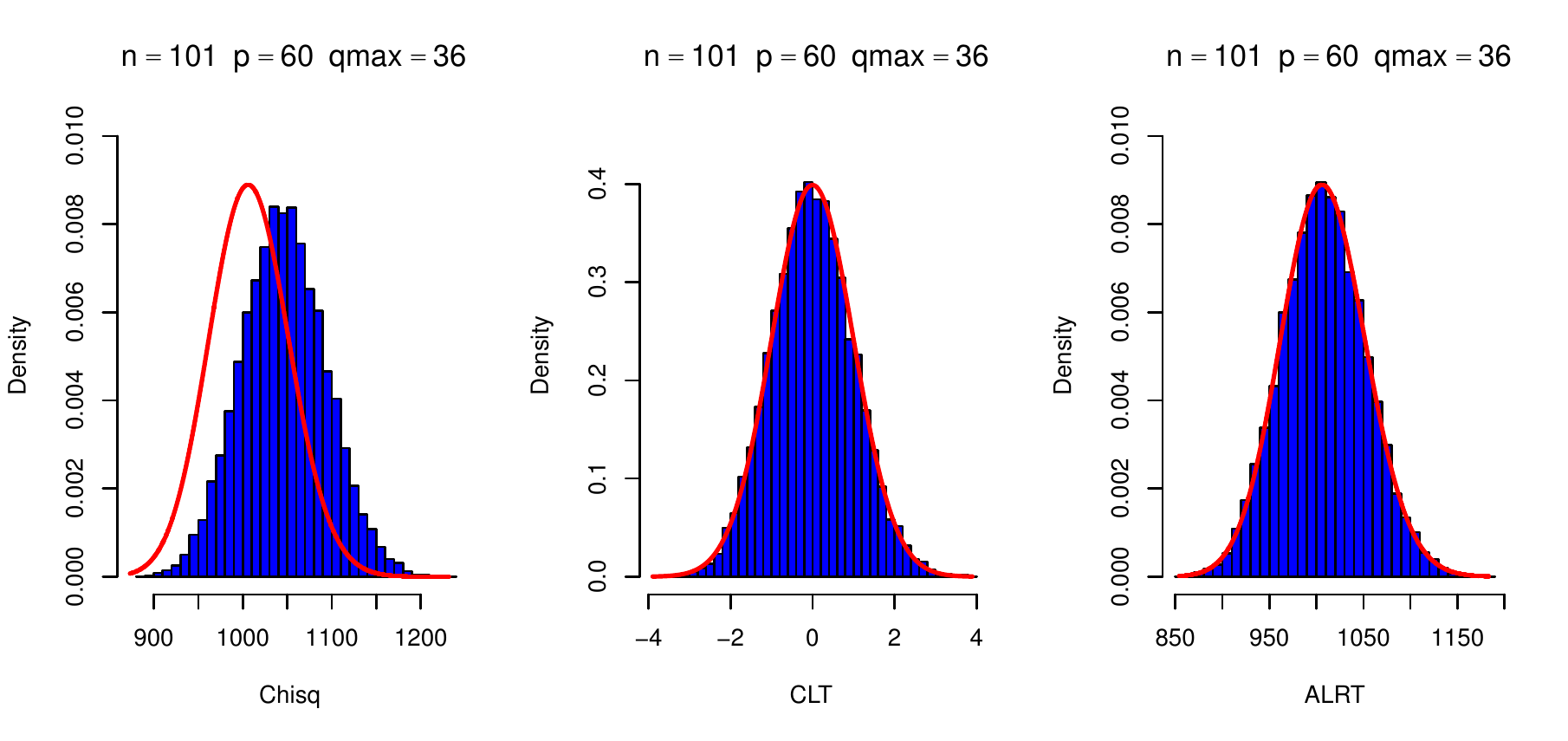}

\includegraphics[height=.2\textheight, width=0.9\textwidth]{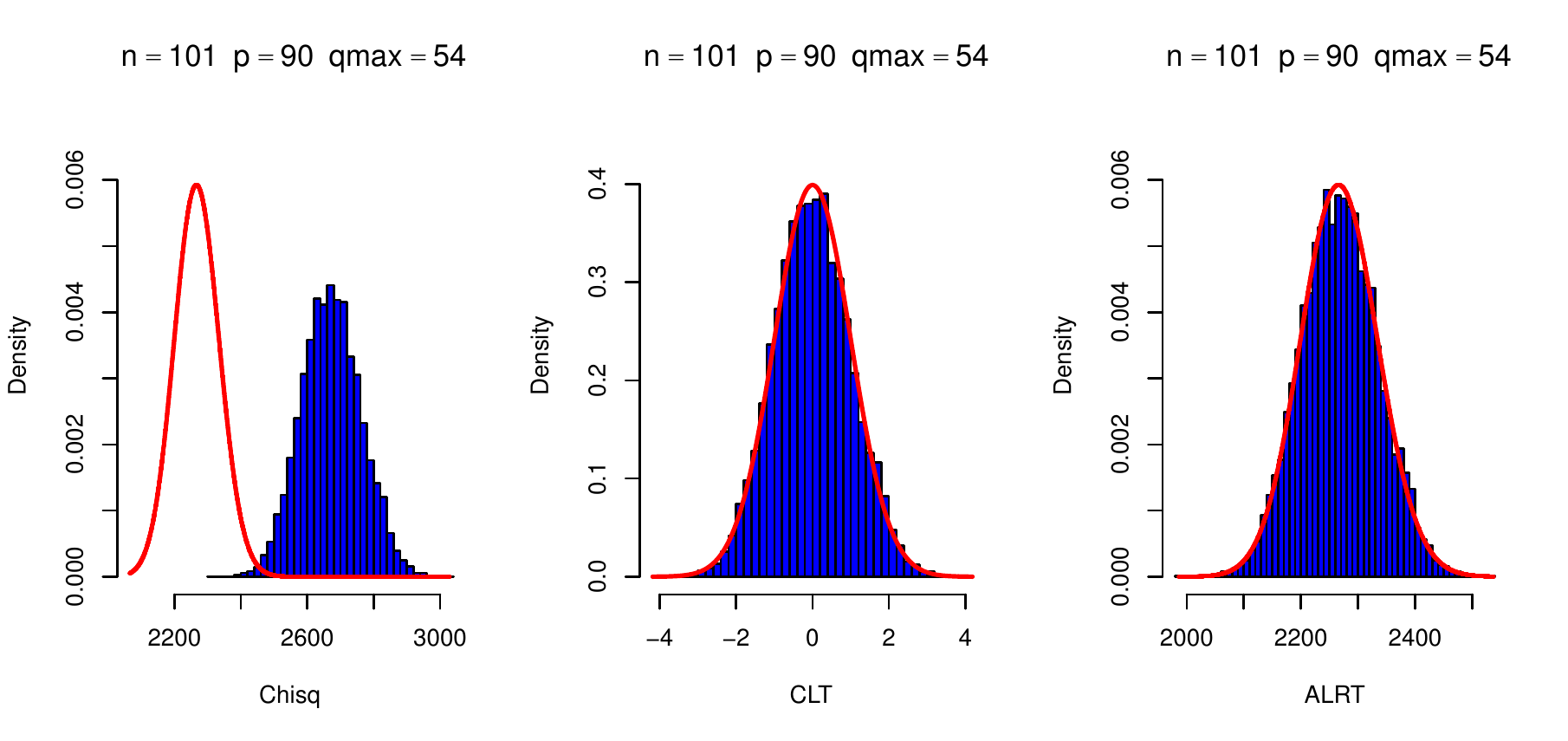}
\caption {Plots of histograms and theoretical density curves (smooth
curves in graphs) for three test statistics based on the classical
chi-square approximation \eqref{classic}, the normal approximation
\eqref{T0}, and the adjusted chi-square approximation in
\eqref{chisquare} with selection of $k=3$ and ratio
$q_1:q_2:q_3=3:1:1$.
 }
 \label{figure1}
\end{figure}

\noindent\textbf{Case a}. We set $k=3$, $p=10, 30, 60$ and $90$ and
use the ratio $q_1:q_2:q_3=3:1:1$.  Figure~\ref{figure1} contains 12
plots in an array with four rows and three columns, and each row
corresponds to one value of $p=10, 30, 60$, and $90$.

\begin{figure}[p!]
\centering
\includegraphics[height=.2\textheight, width=0.9\textwidth]{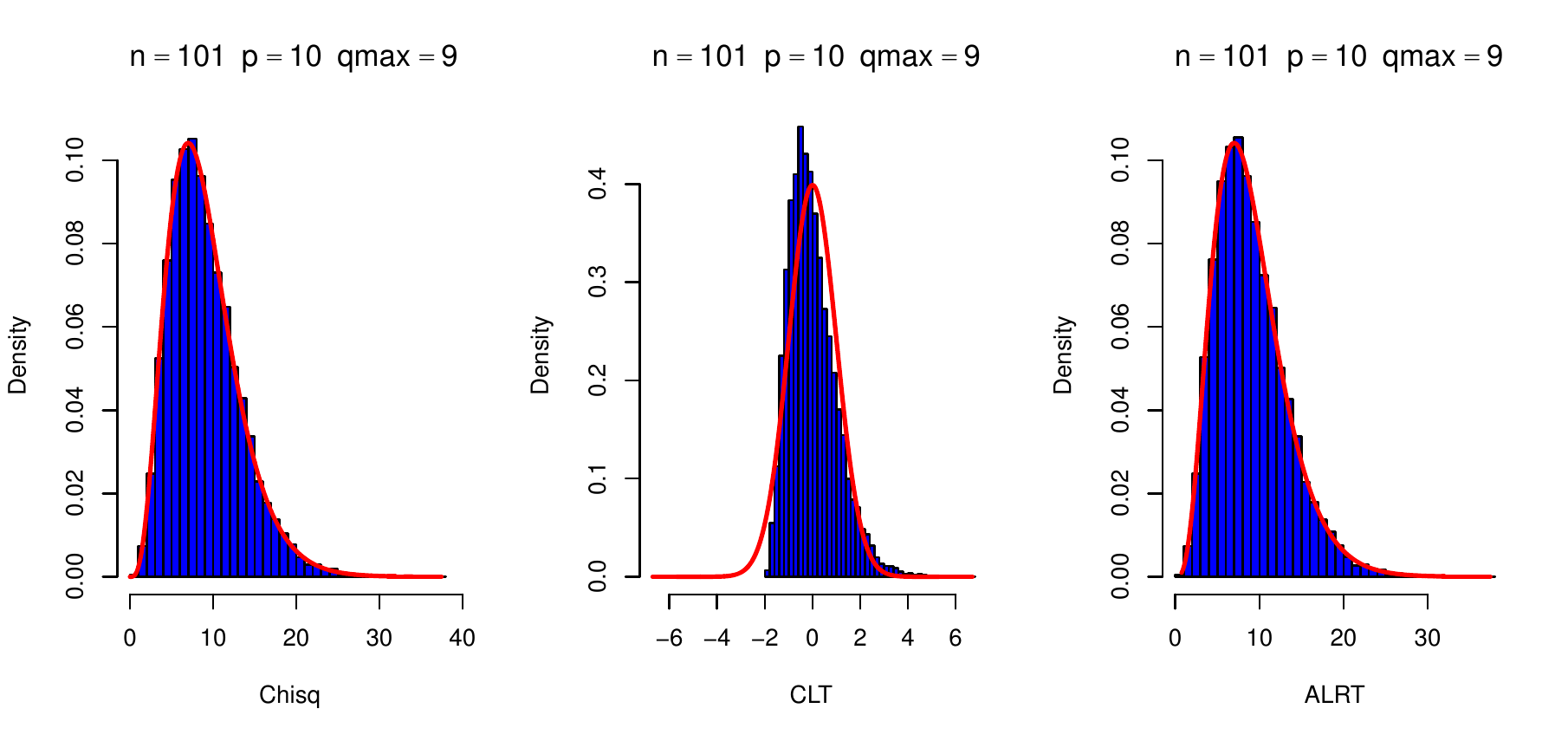}

\includegraphics[height=.2\textheight, width=0.9\textwidth]{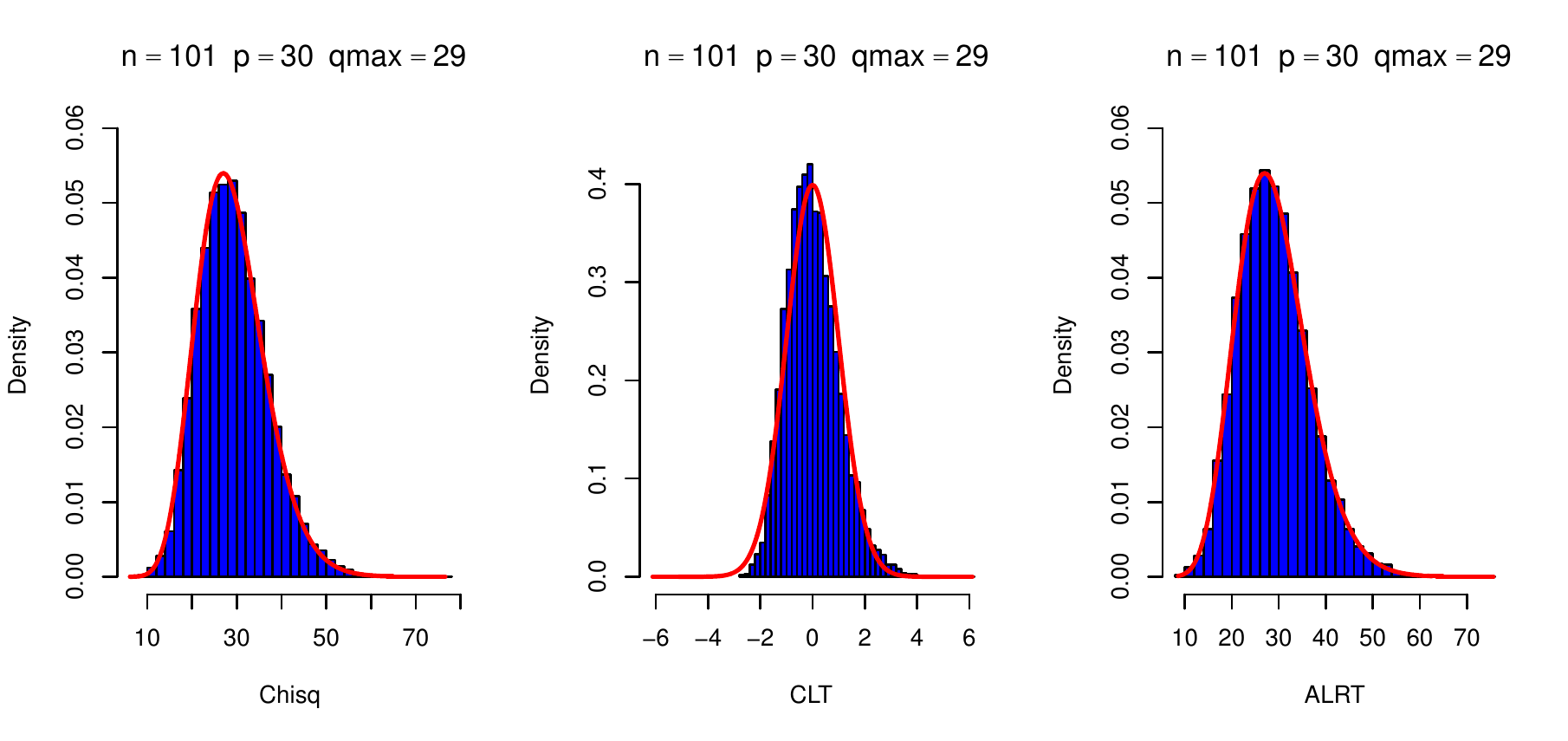}

\includegraphics[height=.2\textheight, width=0.9\textwidth]{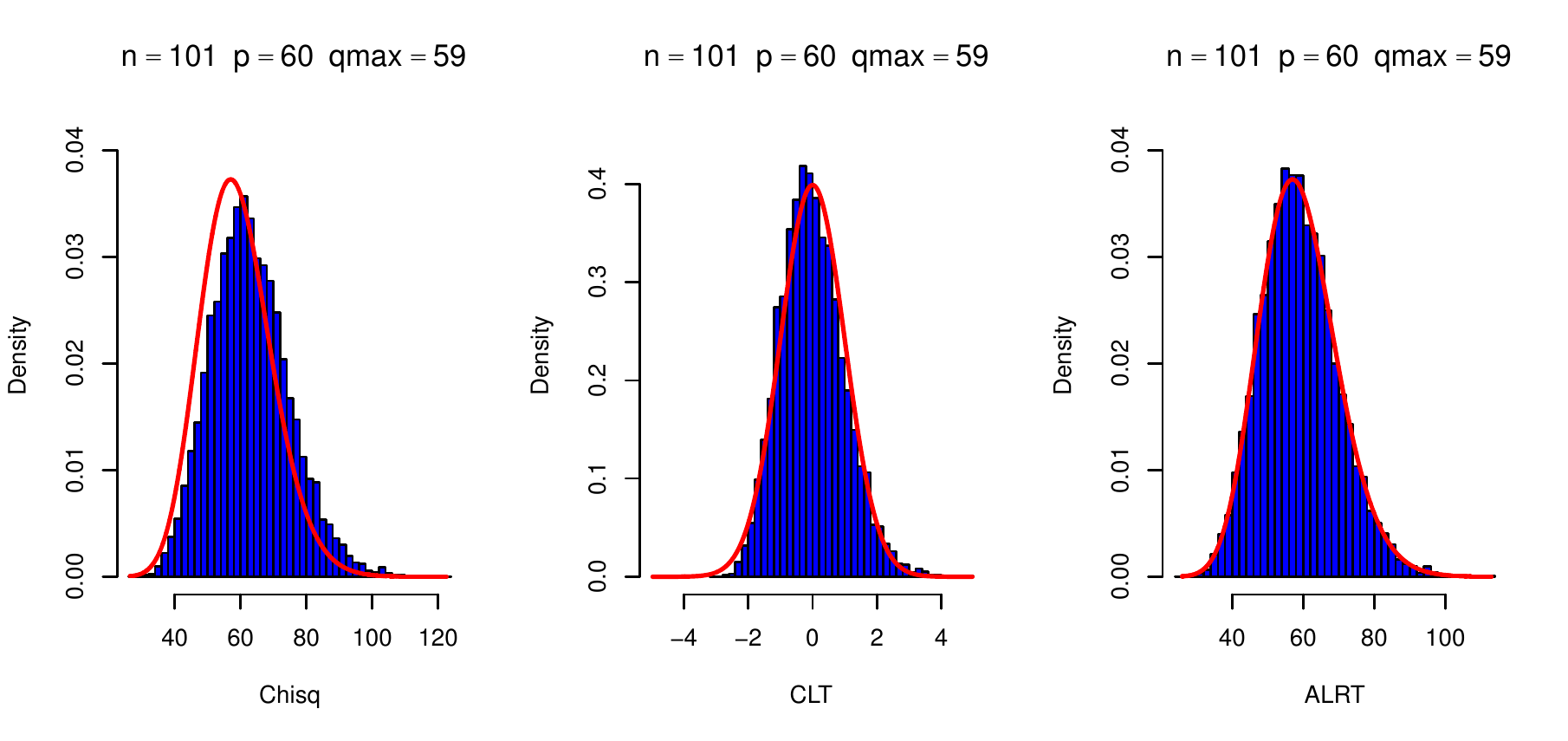}

\includegraphics[height=.2\textheight, width=0.9\textwidth]{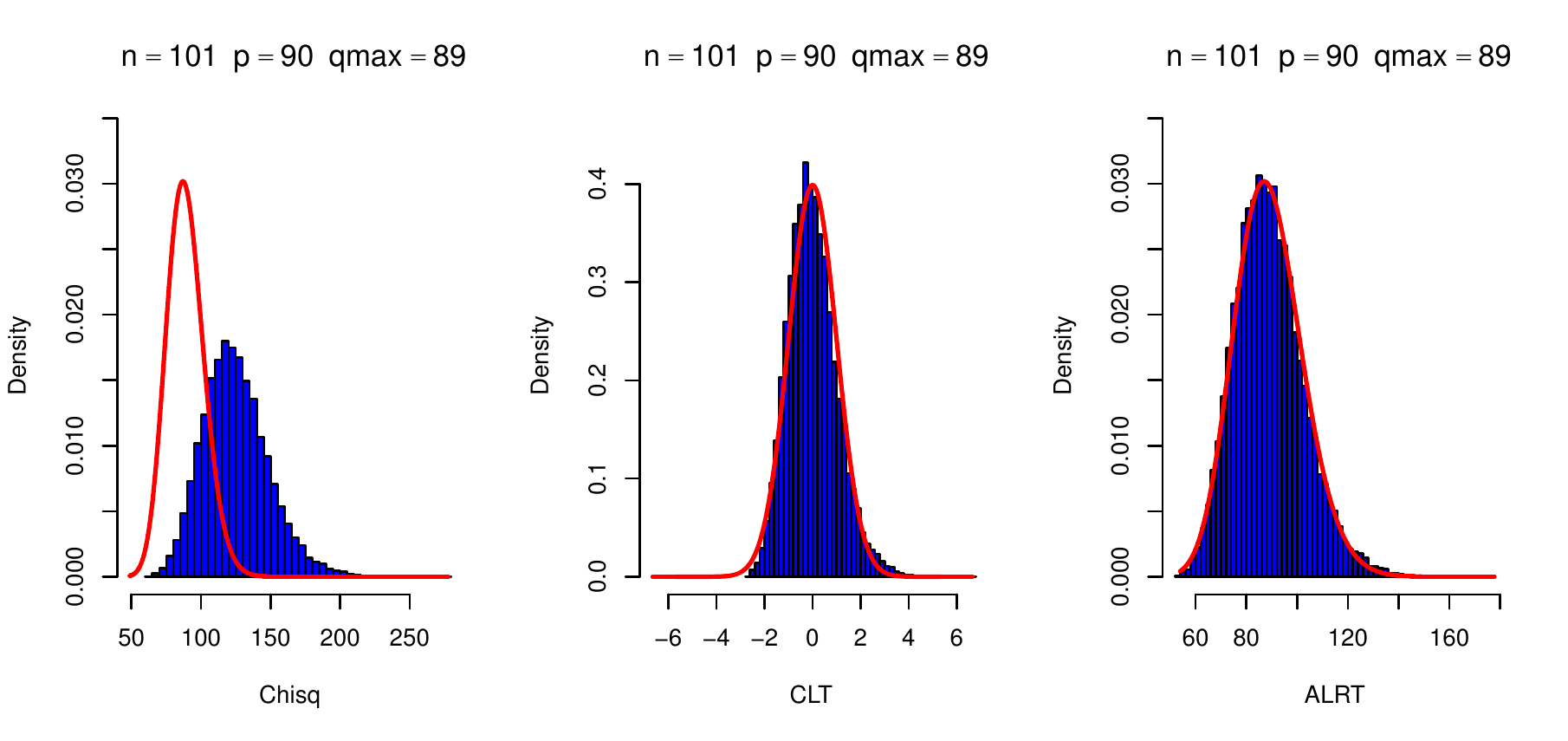}
\caption {Plots of histograms and theoretical density curves (smooth
curves in graphs) for three test statistics based on the classical
chi-square approximation \eqref{classic}, the normal approximation
\eqref{T0}, and the adjusted chi-square approximation in
\eqref{chisquare} with selection of $k=2$, $(q_1, q_2)=(p-1, 1)$ and
$p=10, 30, 60$ and $90$.
 }
 \label{figure2}
\end{figure}

\noindent\textbf{Case b}. We set $k=2$, $q_1=p-1$ and $q_2=1$ for
$p=10, 30, 60$, and $90$, respectively. Figure~\ref{figure2}
contains 12 plots in an array with four rows and three columns, and
each row corresponds to one value of $p=10, 30, 60$, and $90$.

\begin{figure}[p!]
\centering
\includegraphics[height=.2\textheight, width=0.9\textwidth]{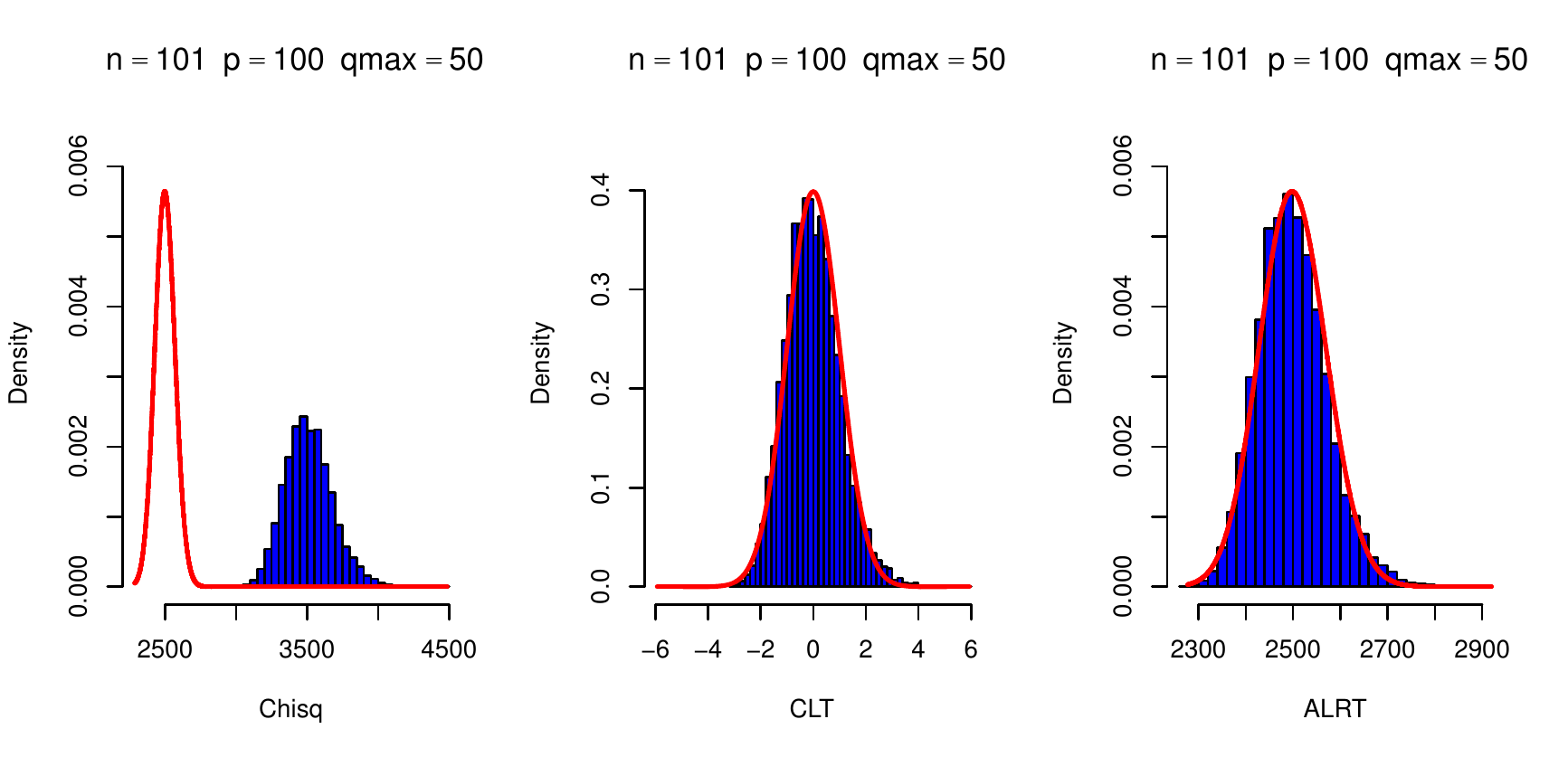}

\includegraphics[height=.2\textheight, width=0.9\textwidth]{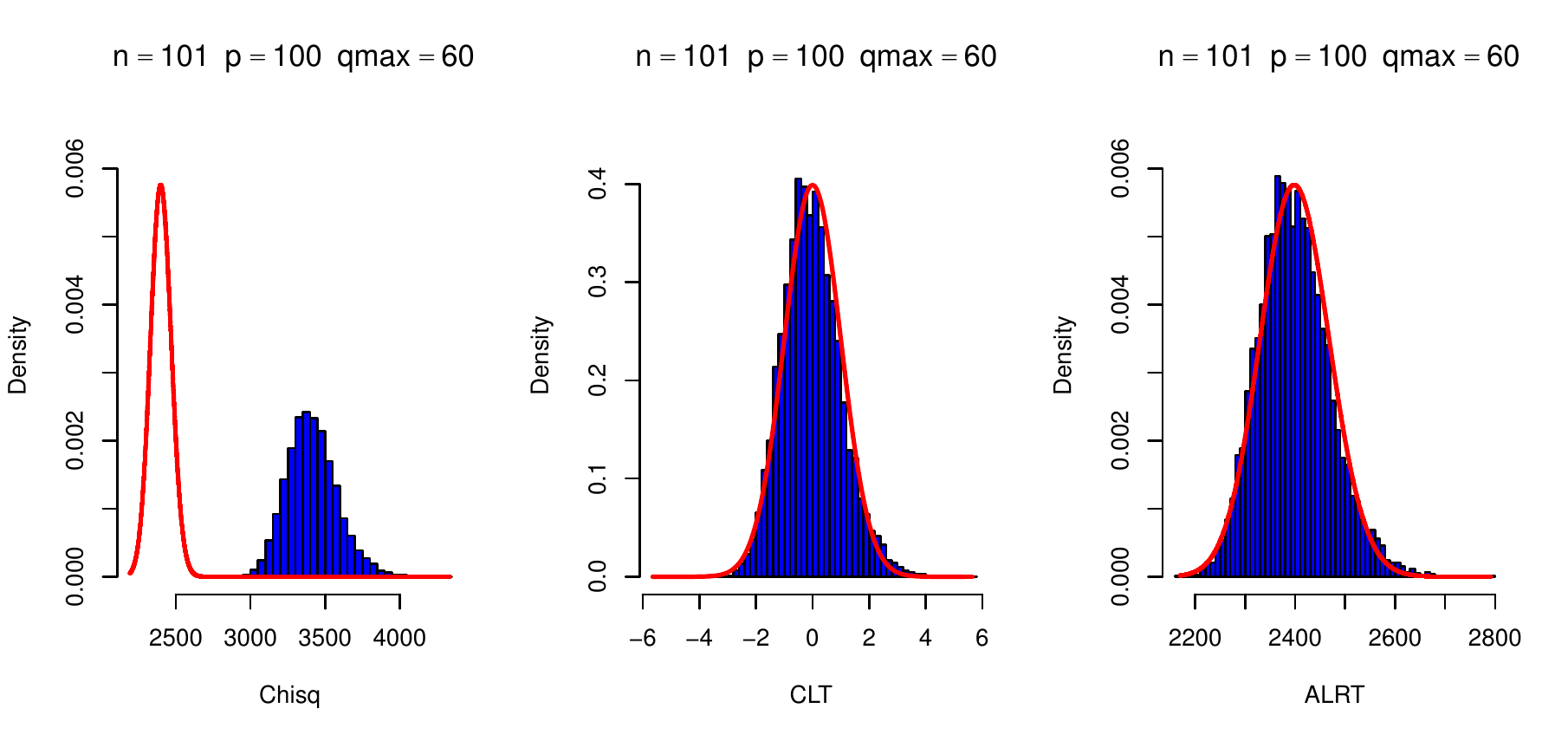}

\includegraphics[height=.2\textheight, width=0.9\textwidth]{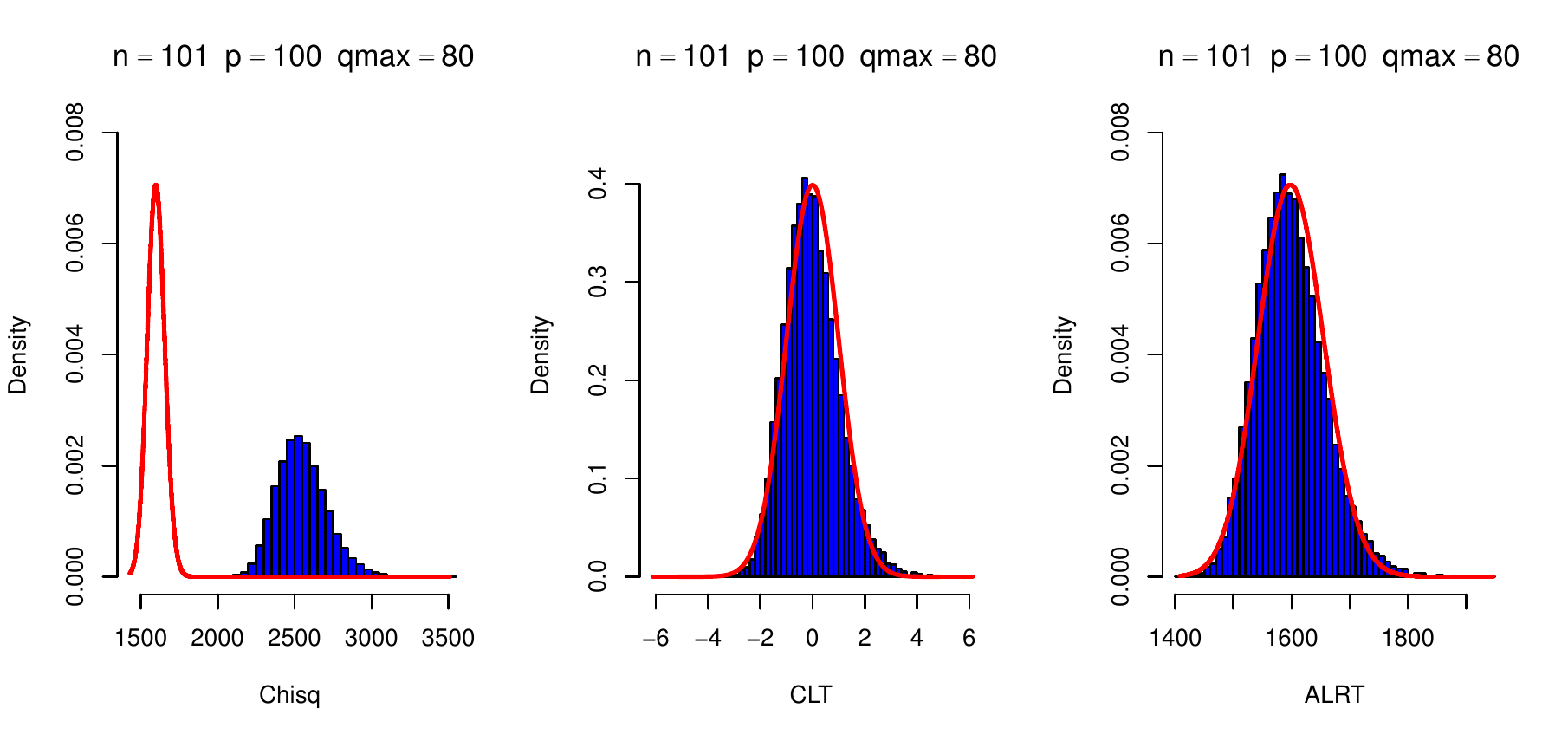}

\includegraphics[height=.2\textheight, width=0.9\textwidth]{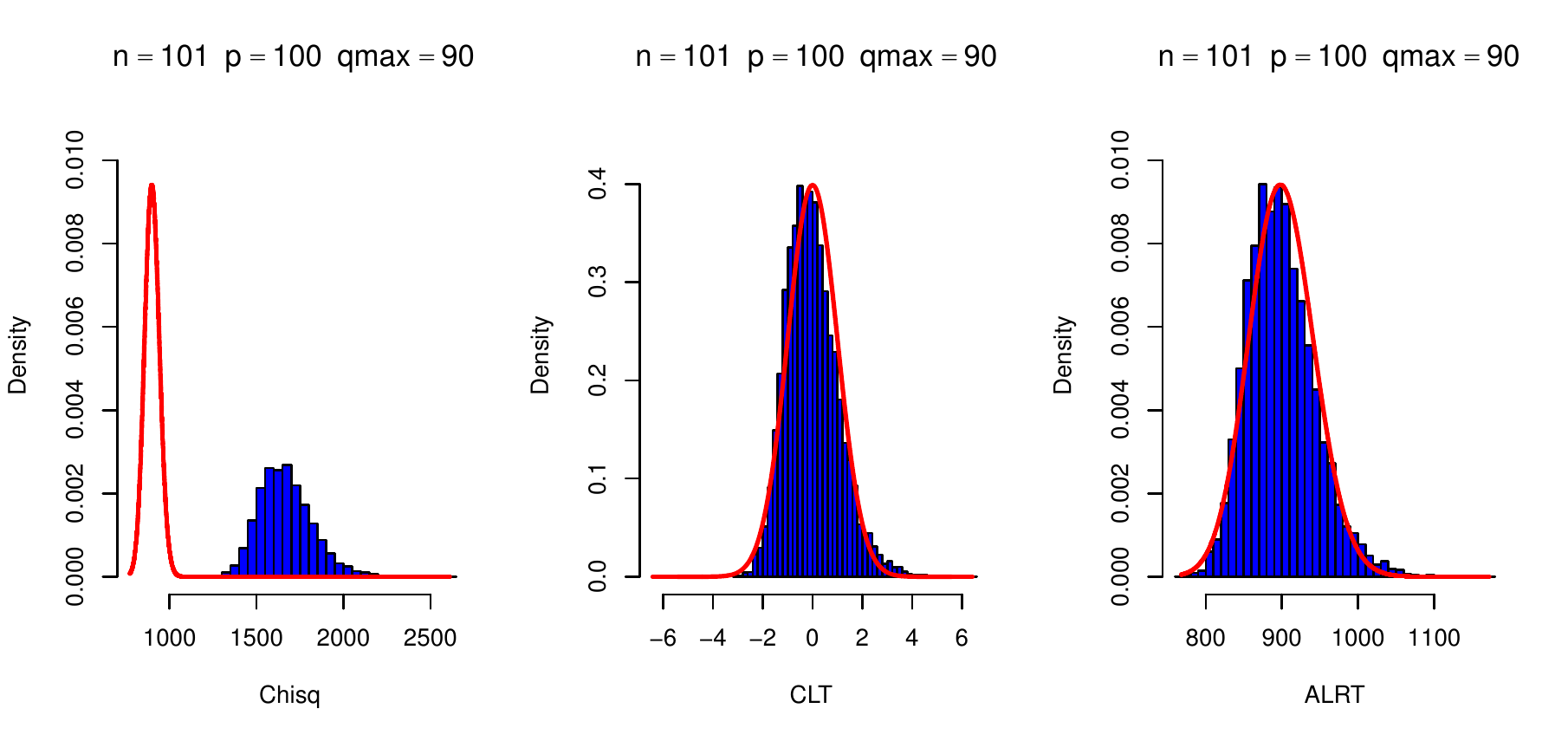}
\caption {Plots of histograms and theoretical density curves (smooth
curves in graphs) for three test statistics based on the classical
chi-square approximation \eqref{classic}, the normal approximation
\eqref{T0}, and the adjusted chi-square approximation in
\eqref{chisquare} with selection of $k=2$, $p=100$ and $(q_1,
q_2)=(50,50)$, $(60,40)$, $(80,20)$, and $(90, 10)$.
 }
 \label{figure3}
\end{figure}

\noindent\textbf{Case c}. We set $k=2$, $p=100$ and choose $(q_1,
q_2)=(50,50)$, $(60, 40)$, $(80, 20)$, and $(90, 10)$, respectively.
Figure~\ref{figure3} contains 12 plots in an array with four rows
and three columns, and each row corresponds to one combination of
$(q_1, q_2)$ with $q_{\max}=50, 60, 80$, and $90$, respectively.

\noindent\textbf{Case d}. We set $n = 101$, $k = 2$, $r =
p-q_{\max}=1$ and and $v = n- p = 1$,  $3$,  $5$ and $10$.

For the first three cases above, we select parameters to maintain a
reasonably large value for $n-q_{\max}$ so that we compare the
performance of classical chi-square approximation \eqref{classic},
the normal approximation \eqref{T0}, and the adjusted chi-square
approximation in \eqref{chisquare}. Under Case \textbf{a}, the
values of $q_i$'s are proportional. Under Cases \textbf{b} and
\textbf{c}, two extreme situations, either $q_{\max}=p-1$ or
$p=n-1$, are considered.

Case \textbf{d} is used to compare the performance of the normal
approximation \eqref{T0}, the adjusted chi-square approximation in
\eqref{chisquare} and the non-normal approximation in
\eqref{non-normal}.  Since we take  $r=1$,  the limit on the
right-hand side of \eqref{non-normal} has only one term, that is,
the limit is $-\log Y_v$, where $Y_v$ is a random variable having a
chi-square distribution with $v$ degrees of freedom.

\begin{figure}[p!]
\centering
\includegraphics[height=.2\textheight, width=0.9\textwidth]{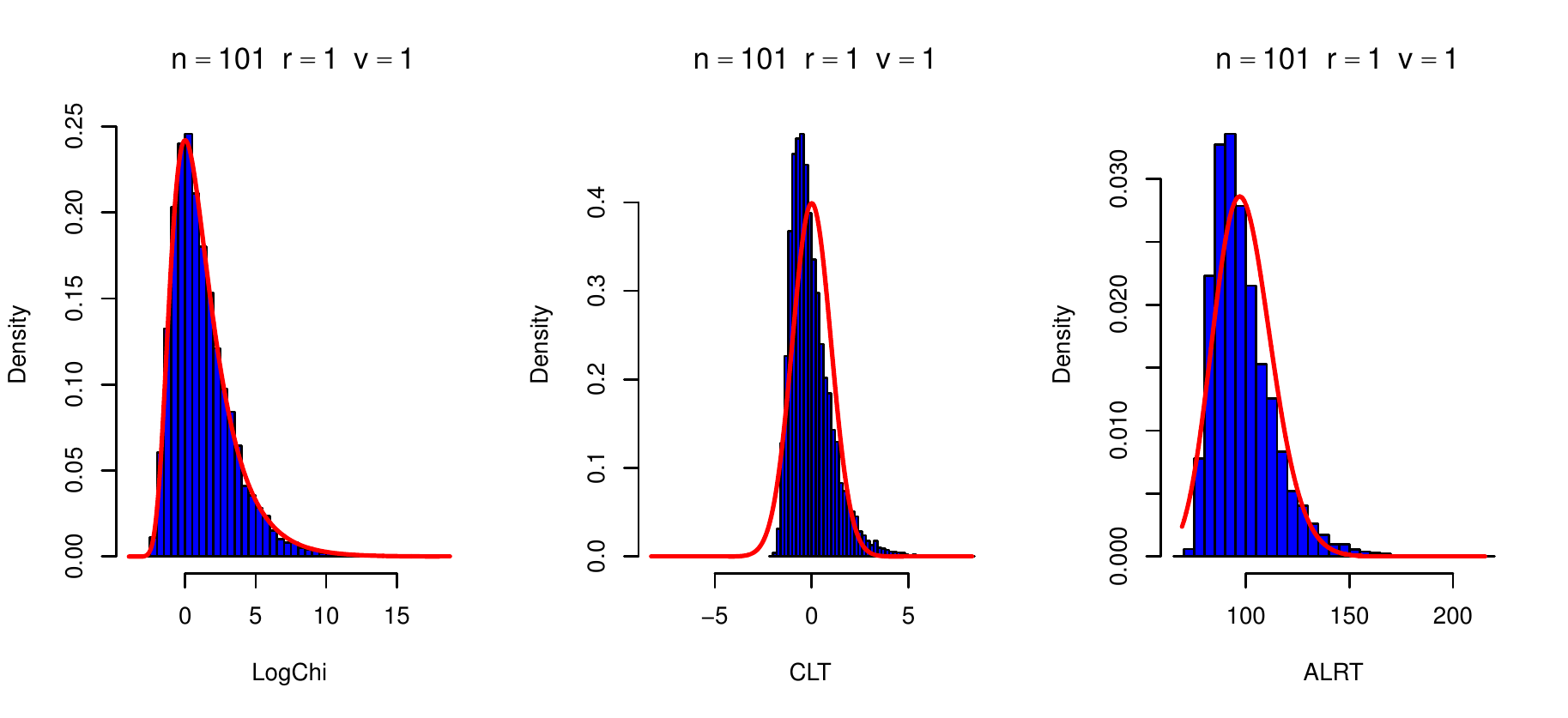}

\includegraphics[height=.2\textheight, width=0.9\textwidth]{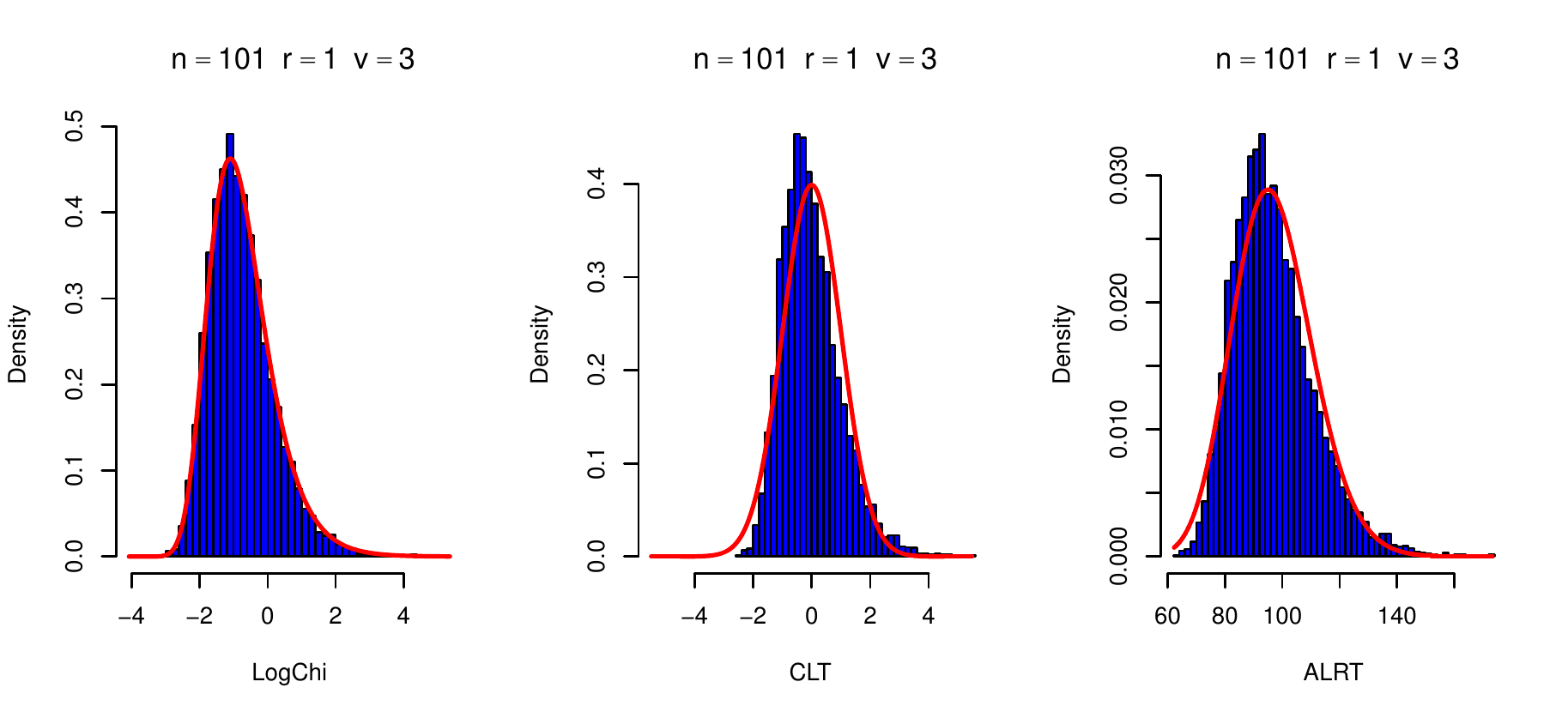}

\includegraphics[height=.2\textheight, width=0.9\textwidth]{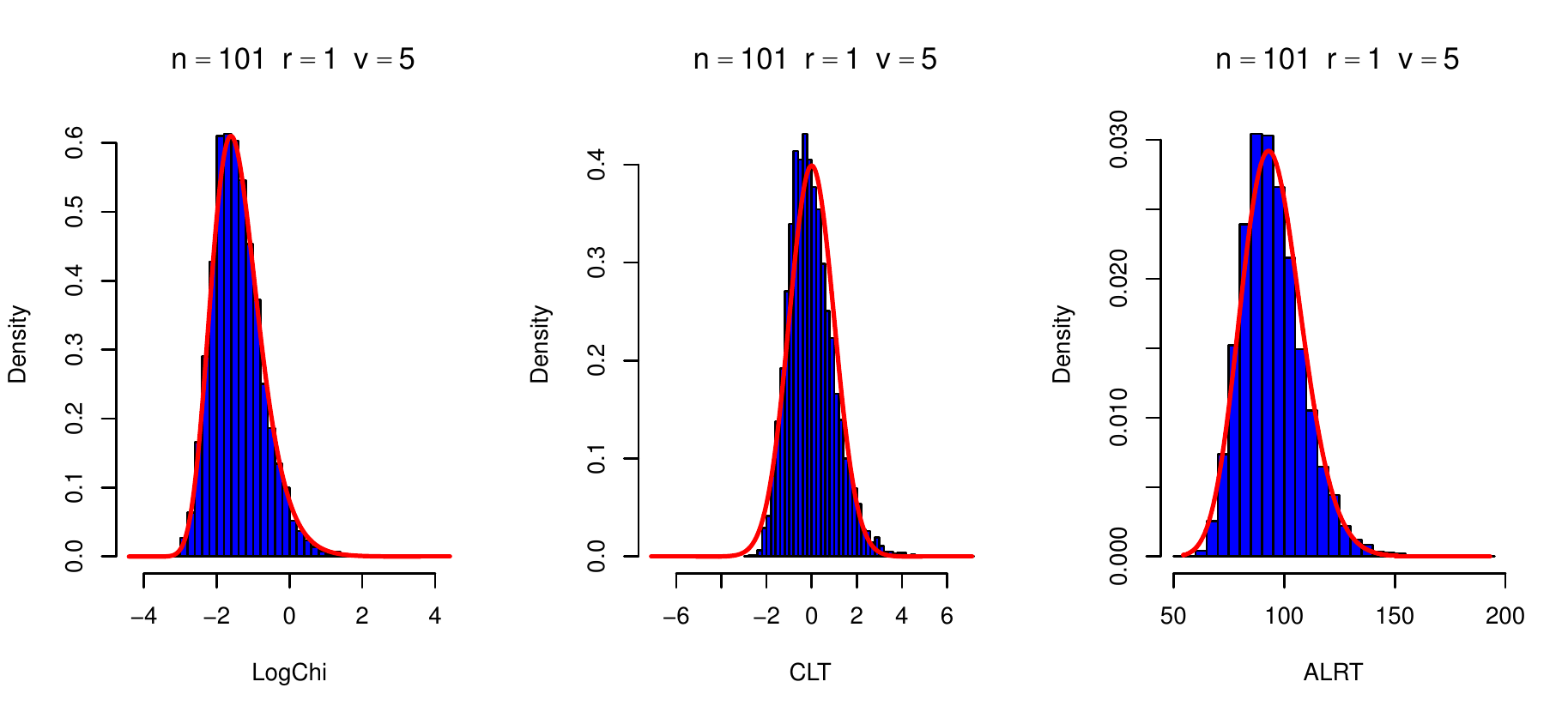}

\includegraphics[height=.2\textheight, width=0.9\textwidth]{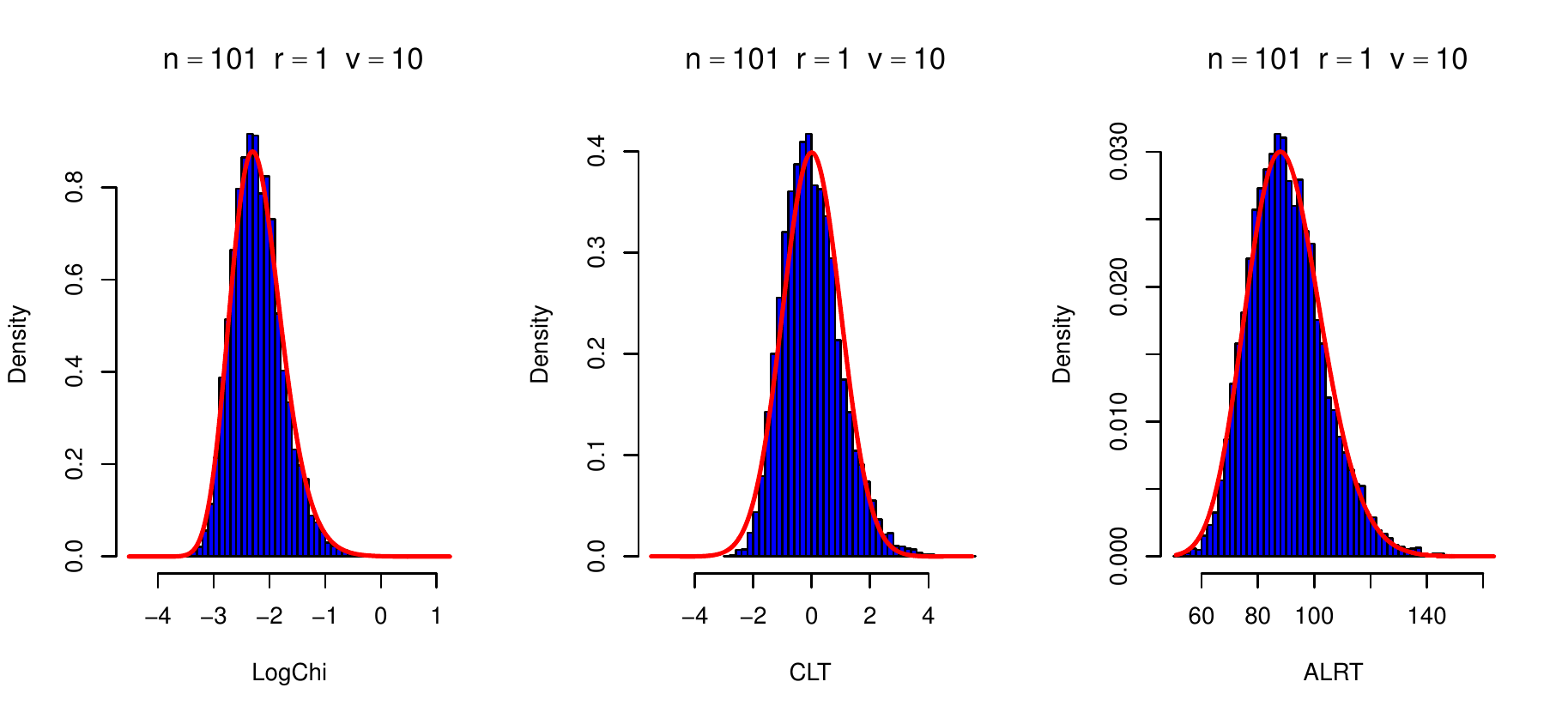}
\caption {Plots of histograms and theoretical density curves (smooth
curves in graphs) for three test statistics based on non-normal
approximation \eqref{non-normal}, the normal approximation
\eqref{T0}, and the adjusted chi-square approximation in
\eqref{chisquare} with selection of $n=101$, $k=2$,
$r=p-q_{\max}=1$, and $v=n-p=1, 3, 5$ and $10$. For the plots in the
first column above, LogChi denotes the non-normal approximation in
\eqref{non-normal}. }
 \label{figure4}
\end{figure}

The results under Cases \textbf{a} to \textbf{d} are given in
Figures~\ref{figure1} to \ref{figure4}, respectively.

Now we summarize our findings from Figures~\ref{figure1} to
\ref{figure4}.

\begin{enumerate}
  \item The classical chi-square approximation (Chisq) works very well for small $p$,
  but it becomes worse with the increase of $p$ and finally departs from the histograms of the test statistic.

  \item When $p$ is small such as $p=10$, the normal approximation (CLT) shows
 lack of fit to the histograms and it is getting better with the increase of $p$.

  \item The adjusted likelihood ratio method (ALRT) works very well for all cases,  that is, for small $p$,
the ALRT behaves like the classical chi-square approximation, while
for large $p$, it performs very well too like the normal
approximation.

  \item In Figure~\ref{figure3}, we select $p=n-1$.  Both the normal approximation
and the adjusted likelihood method show a little bit departure from
 the histograms since $p$ is too close to $n$ and sample size
$n=101$ is not a large sample size. In this case, both the normal
approximation and the adjusted likelihood method can improve when
$n$ is getting larger.

 \item From Figure~\ref{figure4},  when $v$ is small,  the
 non-normal approximation works much better than the normal
 approximation and the adjusted chi-square approximation. When the value of  $v$
  increases from $5$ to $10$, both the normal
 approximation and the adjusted chi-square approximation improve
 significantly.  This implies that one can use the normal
 approximation or the adjusted chi-square approximation when
 $r+v=n-p_{\max}$ is not too small. We note that the exact distribution of the
 limit on the right-hand side of \eqref{non-normal} is not easy to
 obtain if $r=p-q_{\max}>1$.

\end{enumerate}

\clearpage


\subsection{Comparison of adjusted log-LRT and other methods}\label{sub2}

In this subsection, we compare our adjusted log-likelihood ratio
test statistic, i.e., $Z_n$ in \eqref{zn} with other three test
statistics, including two trace criterion test statistics by Jiang
et al.~\cite{JBZ13} and Li et al.~\cite{Li2017} and Schott type
statistics by Bao et al.~\cite{Bao2017}.  We notice that Jiang et
al.~\cite{JBZ13} and Bao et al.~\cite{Bao2017} investigate their
test statistics for test \eqref{t1} for any fixed $k\ge 2$ while Li
et al.~\cite{Li2017} consider test \eqref{t1} for $k=2$ only.  Given
that Li et al.~\cite{Li2017}'s test statistics are suitable for
$k=2$ only, we set $k=2$ for the comparison in this section.

In Jiang et al.~\cite{JBZ13}, a large-dimensional trace criterion
test statistic is defined as
 \[
 L_n=\textbf{tr}(A_{21}A_{11}^{-1}A_{12}A_{22}^{-1}),
 \]
 where $\textbf{tr}(A)$ denotes the trace of matrix $A$. Under conditions $r_{n1}:=q_2/q_1\to r_1\in (0, \infty)$,
 $r_{n2}:=q_2/(n-1-q_2)\to r_2\in (0,\infty)$, and $q_2<n$, it is
 shown that
\begin{equation}\label{T2}
T_2:=\frac{L_n-a_n}{\sqrt{b_n}}\overset{d}\to N(0,1) \mbox{ as } n\to\infty
\end{equation}
under the null hypothesis in \eqref{t1}, where
\[
b_n=\frac{2h_n^2r_{n1}^2r_{n2}^2}{(r_{n1}+r_{n2})^2}, ~a_n=\frac{q_2r_{n2}}{r_{n1}+r_{n2}}, ~h_n=\sqrt{r_{n1}+r_{n2}-r_{n1}r_{n2}}.
\]
Some calculation shows that $a_n=q_1q_2/(n-1)$ and
$b_n=2q_1q_2(n-1-q_1)(n-1-q_2)/(n-1)^4$.

When $k=2$, the test statistic in Bao et al.~\cite{Bao2017} is equal
to
\[
\textbf{tr}(A_{22}^{-1/2}A_{21}A_{11}^{-1}A_{12}A_{22}^{-1/2})=\textbf{tr}(A_{21}A_{11}^{-1}A_{12}A_{22}^{-1}),
\]
which is the same $L_n$. Meanwhile, Theorem 3.1 in \cite{Bao2017}
also implies \eqref{T2}.

To introduce Li et al.~\cite{Li2017}'s test statistic, set for $i,
j=1,2$
\[
\gamma_{ij}=\frac{1}{(n-2)(n+1)}\Big(\textbf{tr}(A_{ij}A_{ji})-\frac{1}{n-1}\textbf{tr}(A_{ii})\textbf{tr}(A_{jj})\Big).
\]
The trace criterion test statistic by \cite{Li2017} is defined as
$\gamma_{12}$. Under the null hypothesis in \eqref{t1} for $k=2$, it
is shown in Li et al.~\cite{Li2017} that
\begin{equation}\label{T3}
T_3:=\sqrt{\frac{(n-2)(n+1)}{2}}\frac{\gamma_{12}}{\sqrt{\gamma_{11}\gamma_{22}}}\overset{d}\to N(0,1) \mbox{ as } n\to\infty
\end{equation}
given that $p=q_1+q_2\to\infty$ as $n\to\infty$ and
\[
0<\lim_{n\to\infty}\frac{1}{p}\textbf{tr}(\mathbf{\Sigma^i})<\infty ~\mbox{ for } i=1, 2, 4.
\]

We notice that the test statistic $T_3$ works only for the case
$k=2$, but condition $p=q_1+q_2\to\infty$ is less restrictive than
those required for other statistics we just discussed.

For any fixed size $\alpha\in (0,1)$, test $T_2$ (or $T_3$) rejects
the null hypothesis in \eqref{t1} if $T_2>z_{\alpha}$ (or
$T_3>z_{\alpha}$), where $z_{\alpha}$ denotes the $\alpha$-level
critical value of the standard normal distribution.

To compare our adjusted log-likelihood ratio test statistic $Z_n$
and test statistics $T_2$ and $T_3$, we assume $q_1>q_2\ge 1$ and
$p=q_1+q_2<n$. Our samples are generated from the populations
similar to those in Jiang et al.~\cite{JBZ13} and Qi et
al.~\cite{QWZ19}. Let $\mathbf{z}=(z_1, \cdots, z_p)'$ be a random
vector whose components are independent normal random variables with
mean $0$ and variance $1$.

\noindent{\bf Model 1.}  $\mathbf{x}=(x_1, \cdots, x_p)'$, where $x_i=(1+c)z_i$ for $i=1,\cdots,p_1$,
$x_{p_1+j}=z_{p_1+j}+cz_j$ for $j=1,\cdots,p_2$, and $c$ is a constant;

\noindent{\bf Model 2.}  $\mathbf{x}=(x_1, \cdots, x_p)'$, where $x_i=(1+c)z_i$ for $i=1,\cdots,p_1$,
$x_{p_1+j}=z_{p_1+j}+cz_j$ for $j=1,\cdots,p_2-1$, $x_p=p^{-1/4}z_p$, and $c$ is a constant.

With different selections of $(q_1, q_2, n, c)$, we draw $10000$
random samples of size $n$ from each model above (Model 1 and Model
2) and then we calculate the empirical sizes of the tests (when
$c=0$) or the powers of the tests (when $c\ne 0$). The nominal level
of size $\alpha$ (type I error) is set to be $0.05$ in the
simulation.

Tables~\ref{model 1} and \ref{model 2} present results for the
numerical comparisons on the three test statistics. From the two
tables, the adjusted log-likelihood ratio test statistic $Z_n$ is
constantly accurate in terms of type I error, and $T_3$ has larger
type I errors than the nominal level $0.05$ for small values of $p$.
In view of empirical powers,  $Z_n$ and $T_2$ are comparable in most
cases while $T_3$ is better than both $Z_n$ and $T_2$ under Model 1.
Under Model 2, $Z_n$ has a slightly larger power than $T_2$ in most
cases and both are significantly better than $T_3$.

\begin{table}[t!]
\caption{Comparisons on Size and Power under Model 1. The sizes and
powers are estimated based on $10,000$ simulations under Model 1,
and the nominal type I errors for all tests are set to be $0.05$.}
\begin{tabular}{|cr|ccc|ccc|ccc|}
  \hline
  & & \multicolumn{3}{|c|}{Size ($c=0$)} & \multicolumn{3}{|c|}{Power ($c=0.1$)} & \multicolumn{3}{|c|}{Power ($c=0.2$)} \\
  \cline{3-11}
  ($q_1,q_2$) & $n$ & $Z_n$ &  $T_2$ & $T_3$ &  $Z_n$&$T_2$ &$T_3$ &  $Z_n$&$T_2$ &$T_3$\\
  \hline
  (6, 4) & 20 & 0.0507& 0.0514& 0.0645& 0.0605& 0.0583& 0.0806& 0.0909& 0.0923& 0.1394\\
         & 50 & 0.0507& 0.0593& 0.0641& 0.0857& 0.0999& 0.1128& 0.2400& 0.2654& 0.3044\\
         &100 & 0.0527& 0.0615& 0.0645& 0.1350& 0.1570& 0.1690& 0.5643& 0.6028& 0.6262\\
  (8, 2) & 20 & 0.0498& 0.0501& 0.0609& 0.0532& 0.0547& 0.0705& 0.0717& 0.0704& 0.1029\\
         & 50 & 0.0487& 0.0581& 0.0641& 0.0699& 0.0838& 0.0927& 0.1499& 0.1719& 0.1983\\
         &100 & 0.0530& 0.0665& 0.0646& 0.0957& 0.1154& 0.1251& 0.3086& 0.3461& 0.3709\\
          \hline
  (18, 12)& 50& 0.0499& 0.0542& 0.0545& 0.1897& 0.1983& 0.2255& 0.8617& 0.8732& 0.9147\\
          &100& 0.0506& 0.0528& 0.0554& 0.1233& 0.1289& 0.1482& 0.5492& 0.5702& 0.6757\\
          &150& 0.0499& 0.0542& 0.0545& 0.1897& 0.1983& 0.2255& 0.8617& 0.8732& 0.9147\\
  (24,  6)& 50& 0.0499& 0.0472& 0.0562& 0.0647& 0.0614& 0.0777& 0.1133& 0.1115& 0.1811\\
          &100& 0.0478& 0.0519& 0.0540& 0.0854& 0.0881& 0.1096& 0.2990& 0.3128& 0.3855\\
          &150& 0.0547& 0.0595& 0.0588& 0.1214& 0.1334& 0.1482& 0.5393& 0.5587& 0.6210\\
           \hline
  (36, 24)&100& 0.0511& 0.0543& 0.0519& 0.1053& 0.1112& 0.1495& 0.4243& 0.4599& 0.6915\\
          &150& 0.0493& 0.0494& 0.0527& 0.1619& 0.1629& 0.2224& 0.8144& 0.8361& 0.9328\\
          &200& 0.0538& 0.0554& 0.0562& 0.2428& 0.2498& 0.3138& 0.9715& 0.9756& 0.9915\\
  (48, 12)&100& 0.0519& 0.0523& 0.0509& 0.0771& 0.0795& 0.1039& 0.2235& 0.2325& 0.3836\\
          &150& 0.0503& 0.0496& 0.0568& 0.1115& 0.1137& 0.1368& 0.4657& 0.4793& 0.6399\\
          &200& 0.0473& 0.0503& 0.0513& 0.1471& 0.1493& 0.1817& 0.7222& 0.7334& 0.8261\\
           \hline
  (60, 40)&150& 0.0499& 0.0495& 0.0502& 0.1248& 0.1362& 0.2202& 0.6453& 0.7046& 0.9407\\
          &200& 0.0524& 0.0552& 0.0555& 0.2036& 0.2140& 0.3101& 0.9358& 0.9497& 0.9940\\
          &300& 0.0532& 0.0547& 0.0526& 0.4046& 0.4153& 0.5247& 1.0000& 1.0000& 1.0000\\
  (80, 20)&150& 0.0501& 0.0514& 0.0553& 0.0907& 0.0959& 0.1343& 0.3269& 0.3543& 0.6420\\
          &200& 0.0521& 0.0529& 0.0527& 0.1287& 0.1325& 0.1795& 0.6093& 0.6287& 0.8362\\
          &300& 0.0507& 0.0519& 0.0502& 0.2162& 0.2203& 0.2791& 0.9453& 0.9500& 0.9844\\
  \hline
\end{tabular}
\label{model 1}
\end{table}


\begin{table}[t!]
\caption{Comparisons on Size and Power under Model 2. The sizes and
powers are estimated based on $10,000$ simulations under Model 2,
and the nominal type I errors for all tests are set to be $0.05$.}
\begin{tabular}{|cr|ccc|ccc|ccc|}
  \hline
  & & \multicolumn{3}{|c|}{Size ($c=0$)} & \multicolumn{3}{|c|}{Power ($c=0.1$)} & \multicolumn{3}{|c|}{Power ($c=0.2$)} \\
  \cline{3-11}
  ($q_1,q_2$) & $n$ & $Z_n$ &  $T_2$ & $T_3$ &  $Z_n$&$T_2$ &$T_3$ &  $Z_n$&$T_2$ &$T_3$\\
  \hline
  (6, 4) & 20 & 0.0507& 0.0514& 0.0658& 0.0797& 0.0798& 0.0834& 0.1965& 0.1991& 0.1525\\
         & 50 & 0.0507& 0.0593& 0.0658& 0.1916& 0.2083& 0.1175& 0.7428& 0.7472& 0.3532\\
         &100 & 0.0527& 0.0615& 0.0647& 0.4467& 0.4774& 0.1880& 0.9938& 0.9945& 0.7142\\
  (8, 2) & 20 & 0.0498& 0.0501& 0.0626& 0.0787& 0.0763& 0.0749& 0.1855& 0.1805& 0.1214\\
         & 50 & 0.0487& 0.0581& 0.0695& 0.1966& 0.2178& 0.1057& 0.7311& 0.7443& 0.2567\\
         &100 & 0.0530& 0.0665& 0.0654& 0.4261& 0.4626& 0.1456& 0.9901& 0.9921& 0.5237\\
          \hline
  (18, 12)& 50& 0.0494& 0.0492& 0.0579& 0.1365& 0.1330& 0.1000& 0.5468& 0.4722& 0.3095\\
          &100& 0.0506& 0.0528& 0.0564& 0.4118& 0.3876& 0.1582& 0.9961& 0.9841& 0.7051\\
          &150& 0.0499& 0.0542& 0.0571& 0.7148& 0.6869& 0.2317& 1.0000& 1.0000& 0.9343\\
  (24,  6)& 50& 0.0499& 0.0472& 0.0570& 0.1378& 0.1279& 0.0852& 0.5213& 0.4103& 0.1983\\
          &100& 0.0478& 0.0519& 0.0536& 0.4234& 0.3949& 0.1132& 0.9937& 0.9718& 0.4257\\
          &150& 0.0547& 0.0595& 0.0601& 0.7320& 0.6968& 0.1587& 1.0000& 1.0000& 0.6789\\
           \hline
  (36, 24)&100& 0.0511& 0.0543& 0.0517& 0.2865& 0.2672& 0.1492& 0.9389& 0.8449& 0.7075\\
          &150& 0.0493& 0.0494& 0.0535& 0.6008& 0.5303& 0.2269& 0.9999& 0.9979& 0.9392\\
          &200& 0.0538& 0.0554& 0.0576& 0.8552& 0.7843& 0.3194& 1.0000& 1.0000& 0.9926\\
  (48, 12)&100& 0.0519& 0.0523& 0.0511& 0.2822& 0.2407& 0.1048& 0.9086& 0.7314& 0.4058\\
          &150& 0.0503& 0.0496& 0.0553& 0.6079& 0.5098& 0.1447& 0.9998& 0.9906& 0.6671\\
          &200& 0.0473& 0.0503& 0.0519& 0.8579& 0.7741& 0.1913& 1.0000& 0.9999& 0.8484\\
           \hline
  (60, 40)&150& 0.0499& 0.0495& 0.0493& 0.3966& 0.3446& 0.2209& 0.9915& 0.9530& 0.9460\\
          &200& 0.0524& 0.0552& 0.0560& 0.7092& 0.5988& 0.3147& 1.0000& 0.9993& 0.9962\\
          &300& 0.0532& 0.0547& 0.0535& 0.9768& 0.9346& 0.5361& 1.0000& 1.0000& 1.0000\\
  (80, 20)&150& 0.0501& 0.0514& 0.0564& 0.3863& 0.3046& 0.1403& 0.9766& 0.8280& 0.6582\\
          &200& 0.0521& 0.0529& 0.0528& 0.7005& 0.5473& 0.1865& 1.0000& 0.9924& 0.8517\\
          &300& 0.0507& 0.0519& 0.0510& 0.9791& 0.9093& 0.2885& 1.0000& 1.0000& 0.9886\\
\hline
\end{tabular}
\label{model 2}
\end{table}

\section{Some lemmas}\label{lemma}

The multivariate gamma function, denoted by $\Gamma_p(x)$, is
defined as
\begin{equation}\label{gamma}
\Gamma_p(x):=\pi^{p(p-1)/4}\prod_{i=1}^p\Gamma(x-\frac{1}{2}(i-1))
\end{equation}
with $x>\frac{p-1}{2}$ from Muirhead~\cite{muirhead82}.


For any positive integers $n$ and $p$ with $1<p<n$, let $q_1,
\cdots, q_k$ be $k$ positive integers such that $p=\sum^k_{i=1}q_i$,
where $k\ge 2$ is an integer which may depend on $n$. Denote
$q_{\max}=\max_{1\le i\le k}q_i$.

For any function $g$ defined over $(0,\infty)$,  set
\begin{equation}\label{Psi}
\Psi_{g,\,n,\, p}(x)=\Delta_{g,\,n,\,p}(x)-\sum^k_{i=1}\Delta_{g,\,
n, \,
q_i}(x)=\sum^p_{j=1}g(\frac{n-j}{2}+x)-\sum^k_{i=1}\sum^{q_i}_{j=1}g(\frac{n-j}{2}+x)
\end{equation}
for $x>-\frac{n-p}{2}$, where $\Delta_{g,\, n, \, q_i}$ is defined
in \eqref{Delta}. For brevity, we omit $q_1, \cdots, q_k$ in the
definition of $\Psi_{g,\,n,\, p}(x)$.

Let $g$ be a differentiable function. Both $\Delta_{g,\,n,\, q}$ and
$\Psi_{g,\,n,\,p}$  are linear functionals in $g$ with following
property
\begin{equation}\label{Psi'}
\frac{d}{dx}\Delta_{g,\,n,\,q}(x)=\Delta_{g',\,n,\,q}(x),~~~\frac{d}{dx}\Psi_{g,\,n,\,q}(x)=\Psi_{g',\,
n,\, q}(x).
\end{equation}

Now we set $\beta_{nr}(x)=\Psi_{g,\, n,\, p}(x)$ when $g(x)=1/x^{r}$
for $r\ge 1$, that is, we define
\begin{equation}\label{beta-nr}
\beta_{nr}(x)=\sum^p_{j=1}\frac{1}{\big(\frac{n-j}{2}+x\big)^r}-\sum^k_{i=1}\sum^{q_i}_{j=1}\frac{1}{\big(\frac{n-j}{2}+x\big)^r},
~~x>-\frac{n-p}{2}.
\end{equation}
Note that $\beta_{nr}(x)\ge 0$ since $g(x)=1/x^{r}$, $x>0$ is a
decreasing function over $(0,\infty)$.

Define
\[
s(x)=(x-\frac{1}{2})\log(x)-x, ~x>0
\]
and set
\begin{equation}\label{hx}
h(x)=\log\Gamma(x)-s(x),~~~x>0.
\end{equation}
Then we can verify that
\begin{equation}\label{s'''}
s'(x)=\log
x-\frac{1}{2x},~s''(x)=\frac{1}{x}+\frac{1}{2x^2},~s'''(x)=-\frac{1}{x^2}-\frac{1}{x^3},
\end{equation}
and for some constant $C>0$
\begin{equation}\label{h-bound}
|h''(x)|\le \frac{C}{x^3}, ~~~x\ge \frac14;
\end{equation}
See Lemma 4.4 in Guo and Qi~\cite{GuoQi21}.

For the digamma function $\psi$ defined in \eqref{digamma}, it
follows from Formula 6.3.18
in Abramowitz and
Stegun~\cite{Abramowitz1972} that
\begin{equation}\label{1storder}
\psi(x)=\log x-\frac{1}{2x}-\frac{1}{12x^2}+O(\frac{1}{x^4})
\end{equation}
as $x \rightarrow \infty$.


From now on, we adopt the following notation in our lemmas and our
proofs. For any two sequences, $\{a_n\}$ and $\{b_n\}$ with $b_n>0$,
notation $a_n=o(b_n)$ implies $\lim_{n\to\infty}a_n/b_n=0$, and
notation $a_n=O(b_n)$ means $a_n/b_n$ is uniformly bounded.

We first introduce the formula for the $h$-th moment of $W_n$.

\begin{lemma} \label{lemma1}
(Theorem 11.2.3 in Muirhead~\cite{muirhead82}) Let
$p=\sum_{i=1}^{k}q_k$ and $W_n$ be Wilk's likelihood ratio
statistics defined in \eqref{equ1}. Then, under the null hypothesis
in \eqref{t1}, we have
\begin{equation}\label{MGF}
\mathbf{E}(W_n^h)=\frac{\Gamma_p(\frac{n-1}{2}+h)}{\Gamma_p(\frac{n-1}{2})}\prod_{i=1}^k
\frac{\Gamma_{q_i}(\frac{n-1}{2})}{\Gamma_{q_i}(\frac{n-1}{2}+h)}
\end{equation}
for any $h>\frac{p-n}{2}$, where $\Gamma_p(x)$ is defined in
\eqref{gamma}.
\end{lemma}

Next,  we introduce a distributional representation for $W_n$.

\begin{lemma} \label{lemma1b} (Theorem
11.2.4 in Muirhead~\cite{muirhead82})
 Let $p=\sum_{i=1}^{k}q_k$ and $W_n$ be Wilk's likelihood ratio statistics defined in \eqref{equ1}. Then, under the null hypothesis
in \eqref{t1}, $W_n$ has the same distribution as
\[
\prod_{i=2}^k\prod^{q_i}_{j=1}V_{ij},
\]
where $q_i^*=\sum^{i-1}_{j=1}q_j$ for $2\le i\le k$, the $V_{ij}$'s
are independent random variables and $V_{ij}$ has a
beta($\frac12(n-q_i^*-j)$, $\frac12q_i^*$) distribution.
\end{lemma}

\begin{lemma}\label{sum-q-terms} (Lemma 5 in Qi et al. ~\cite{QWZ19}) As  $n\to\infty$,
\begin{equation}\label{var}
\sum^q_{i=1}(\frac1{n-i}-\frac{1}{n-1})=-\frac{q}{n}-\log(1-\frac{q}{n})+O(\frac{q}{n(n-q)}),
\end{equation}
\begin{equation}\label{mean}
\sum^q_{i=1}\Big(\log(n-i)-\log(n-1)\Big)=(q-n+\frac12)\log(1-\frac{q}{n})-\frac{(n-1)q}{n}+O(\frac{q}{n(n-q)})
\end{equation}
uniformly over $1\le q< n$.
\end{lemma}


\begin{lemma} \label{lemma-beta} With $\beta_{nr}$ defined in
\eqref{beta-nr}, we have
\begin{equation}\label{beta1}
\frac{2p}{3n}\log\big(1+\frac{p-q_{\max}}{3(n-p)}\big)\le
\beta_{n1}(0)\le
\frac{4p}{n}\log\big(1+\frac{2(p-q_{\max})}{n-p}\big).
\end{equation}
\begin{equation}\label{beta2}
\beta_{n2}(x)\le \frac{32\beta_{n1}(0)}{n-p},
\end{equation}
\begin{equation}\label{beta2-bound}
\beta_{n2}(x)\le \frac{8192p(p-q_{\max})}{n(n-p)(n-q_{\max})}
\end{equation}
and
\begin{equation}\label{beta3}
\beta_{n3}(x)\le \frac{192\beta_{n1}(0)}{(n-p)^2}
\end{equation}
for all $|x|\le \frac{n-p}{4}$.
\end{lemma}

\noindent\textit{Proof}. Without loss of generality, assume $q_1\ge
q_2\ge \cdots\ge q_k$. Therefore, $q_1=q_{\max}=\max_{1\le i\le
k}q_i$. Set $q_i^*=\sum^{i-1}_{j=1}q_j$ for $2\le i\le k+1$.  Then
we see that
\begin{eqnarray}\label{beta1=}
\beta_{n1}(0)&=&2\sum_{j=1}^{p}\frac{1}{n-j}-\sum_{i=1}^{k}\sum_{j=1}^{q_i}\frac{1}{n-j}\nonumber\\
&=&2\sum_{i=2}^{k}\sum^{q_i}_{j=1}\big(\frac{1}{n-q_i^*-j}-\frac{1}{n-j}\big)\nonumber\\
&=&2\sum_{i=2}^{k}\sum^{q_i}_{j=1}\frac{q_i^*}{(n-q_i^*-j)(n-j)}.
\end{eqnarray}
Note that $q_i^*\ge q_1=q_{\max}$, and $q_i^*+q_i\le p$ for all
$2\le i\le k$.  Moreover, $q_i\le p/2\le n/2$ for all $2\le i\le k$,
and thus $n>n-j\ge n/2$ for $1\le j\le q_i$, $2\le i\le k$. We have
\begin{eqnarray*}
\beta_{n1}(0)
&\le &\frac{4p}{n}\sum_{i=2}^{k}\sum^{q_i}_{j=1}\frac{1}{n-q_i^*-j}\\
&=&\frac{4p}{n}\sum_{j=1}^{p-q_{\max}}\frac{1}{n-q_{\max}-j}\\
&\le &\frac{4p}{n}\int_{0.5}^{p-q_{\max}+0.5}\frac{1}{n-q_{\max}-x}dx\\
&=&\frac{4p}{n}\log\big(\frac{n-q_{\max}-0.5}{n-p-0.5}\big)\\
&=&\frac{4p}{n}\log\big(1+\frac{2(p-q_{\max})}{n-p}\big),
\end{eqnarray*}
which gives the upper bound in \eqref{beta1}.

To verify the lower bound in \eqref{beta1}, define $m=\min\{i\ge 2:
q_i^*\ge p/3\}$.  We see that $m=2$ if $q_1=q_{\max}\ge p/3$. If
$q_1<p/3$, then $q_k\le \cdots\le q_2\le q_1<p/3$, which implies
that $k\ge 4$, $2<m\le k$, $q_{m-1}^*< p/3$, $p/3\le q_m^*<2p/3$,
and thus, $p-q_m^*\ge p/3\ge (p-q_{\max})/3$. When $q_1<p/3$, it is
trivial that $p-q_m^*\ge p/3\ge (p-q_{\max})/3$. Then we conclude
from \eqref{beta1=} that
\begin{eqnarray*}
\beta_{n1}(0)
&\ge &2\sum_{i=m}^{k}\sum^{q_i}_{j=1}\frac{q_i^*}{(n-q_i^*-j)(n-j)}\\
&\ge &\frac{2p}{3n}\sum_{i=m}^{k}\sum^{q_i}_{j=1}\frac{1}{n-q_i^*-j}\\
&=&\frac{2p}{3n}\sum_{j=1}^{p-q_m^*}\frac{1}{n-q_m^*-j}\\
&\ge &\frac{2p}{3n}\int_{0}^{p-q_m^*}\frac{1}{n-q_m^*-x}dx\\
&=&\frac{2p}{3n}\log\big(\frac{n-q_m^*}{n-p}\big)\\
&=&\frac{2p}{3n}\log\big(1+\frac{p-q_m^*}{n-p}\big)\\
&=&\frac{2p}{3n}\log\big(1+\frac{p-q_{\max}}{3(n-p)}\big),
\end{eqnarray*}
proving \eqref{beta1}.

We can also verify that
\begin{equation}\label{beta2=}
\beta_{n2}(x)=4\sum_{i=2}^{k}\sum^{q_i}_{j=1}\big(\frac{1}{(n-q_i^*-j+2x)^2}-\frac{1}{(n-j+2x)^2}\big)
\end{equation}
and
\begin{equation}\label{beta3=}
\beta_{n3}(x)=8\sum_{i=2}^{k}\sum^{q_i}_{j=1}\big(\frac{1}{(n-q_i^*-j+2x)^3}-\frac{1}{(n-j+2x)^3}\big).
\end{equation}

For any $x$ with $|x|\le (n-p)/4$,  we have for any $1\le j\le q_i$,
$2\le i\le k$
\begin{eqnarray}\label{beta2x}
\frac{1}{(n-q_i^*-j+2x)^2}-\frac{1}{(n-j+2x)^2}
&=&\frac{q_i^*(2n-q_i^*-2j+4x)}{(n-q_i^*-j+2x)^2(n-j+2x)^2}\nonumber\\
&\le&\frac{2q_i^*}{(n-q_i^*-j+2x)^2(n-j+2x)}\nonumber\\
&\le&\frac{16q_i^*}{(n-q_i^*-j)^2(n-j)}\\
&\le&\frac{16}{n-p}\frac{q_i^*}{(n-q_i^*-j)(n-j)},\nonumber
\end{eqnarray}
which, together with \eqref{beta1=}, implies
\[
\beta_{n2}(x)\le
\frac{64}{n-p}\sum_{i=2}^{k}\sum^{q_i}_{j=1}\frac{q_i^*}{(n-q_i^*-j)(n-j)}=\frac{32\beta_{n1}(0)}{n-p},
\]
proving \eqref{beta2}.

Likewise, for any $x$ with $|x|\le (n-p)/4$,  we have for any $1\le
j\le q_i$, $2\le i\le k$
\begin{eqnarray*}
\frac{1}{(n-q_i^*-j+2x)^3}-\frac{1}{(n-j+2x)^3}
&\le&\frac{3q_i^*}{(n-q_i^*-j+2x)^3(n-j+2x)}\\
&\le&\frac{48q_i^*}{(n-q_i^*-j)^3(n-j)}\\
&\le&\frac{48}{(n-p)^2}\frac{q_i^*}{(n-q_i^*-j)(n-j)},
\end{eqnarray*}
which coupled with \eqref{beta3=} yields \eqref{beta3}.

Now we will show \eqref{beta2-bound}. We see that $\max_{2\le i\le
k}q_i\le n/2$ since $\max_{2\le i\le k}q_i+q_1\le p< n$. This
implies $n-j\ge n/2$ for $1\le j\le q_i$, $2\le i\le k$. Then it
follows from \eqref{beta2=} and \eqref{beta2x} that
\[
\beta_{n2}(x)\le
64\sum^k_{i=2}\sum^{q_i}_{j=1}\frac{16q_i^*}{(n-q_i^*-j)^2(n-j)}\le
\frac{2048p}{n}\sum^k_{i=2}\sum^{q_i}_{j=1}\frac{1}{(n-q_i^*-j)^2}=\frac{2048p}{n}\sum^{p}_{j=q_1+1}\frac{1}{(n-j)^2}.
\]
Since the sum in the last expression is dominated by the integral
\[
\int^{p+0.5}_{q_1+0.5}\frac{1}{(n-y)^2}dy\le
\frac{1}{n-p-0.5}-\frac{1}{n-q_1-0.5}=
\frac{p-q_1}{(n-p-0.5)(n-q_1-0.5)}\le
\frac{4(p-q_{\max})}{(n-p)(n-q_{\max})},
\]
we obtain \eqref{beta2-bound}. \eop

\begin{lemma} \label{control} There exists a universal constant
$D$ such that
\begin{equation}\label{control-h}
|\Psi_{h'',\,n,\,p}(x)|\le \frac{Dp}{n(n-p)^2}
\end{equation}
and
\begin{equation}\label{control-s} |\Psi_{s''',\,n,\,p}(x)|\le
\frac{32\beta_{n1}(0)}{n-p}
\end{equation}
uniformly  over $|x|\le \frac{n-p}{4}$ for all large $n$.
\end{lemma}

\noindent\textit{Proof.} For integers $q$ and $n$ with $1\le q<n$,
\[
\sum_{j=1}^{q}\frac{1}{(n-j)^3}\le
\int^{q+0.5}_{0.5}\frac{1}{(n-x)^3}dx=\frac{1}{2}\big(\frac{1}{(n-q-0.5)^2}-\frac{1}{(n-0.5)^2}\big)\le\frac{q}{(n-q-0.5)^2(n-0.5)}\le\frac{8q}{n(n-q)^2}.
\]
We apply the above inequality to $q=q_1, \cdots, q_k$ and $p$. Then
we obtain that
\begin{equation}\label{part1}
\sum_{j=1}^{p}\frac{1}{(n-j)^3}\le \frac{8p}{n(n-q)^2}
\end{equation}
and
\begin{equation}\label{part2}
\sum^k_{i=1}\sum_{j=1}^{q_i}\frac{1}{(n-j)^3}\le
\sum^k_{i=1}\frac{8q_i}{n(n-q_i)^2}\le
\sum^k_{i=1}\frac{8q_i}{n(n-p)^2}=\frac{8p}{n(n-p)^2},
\end{equation}
where we have used the fact that $p=\sum^k_{i=1}q_i$.

In view of \eqref{h-bound}, we have
\begin{eqnarray*}
|\Psi_{h'', n,
p}(x)|&\le&C\sum^p_{j=1}\frac{1}{(\frac{n-j}{2}+x)^3}+C\sum^k_{i=1}\sum^{q_i}_{j=1}\frac{1}{(\frac{n-j}{2}+x)^3}\\
&=&8C\sum^p_{j=1}\frac{1}{(n-j+2x)^3}+8C\sum^k_{i=1}\sum^{q_i}_{j=1}\frac{1}{(n-j+2x)^3}\\
&\le&64C\sum^p_{j=1}\frac{1}{(n-j)^3}+64C\sum^k_{i=1}\sum^{q_i}_{j=1}\frac{1}{(n-j)^3}
\end{eqnarray*}
uniformly over $|x|\le \frac{n-p}{4}$. In the last step above, we
have used the inequality $n-j+2x\ge (n-j)/2$ for all $1\le j\le p$
and $|x|\le (n-p)/4$.  By combining \eqref{part1} and \eqref{part2},
we obtain \eqref{control-h} with $D=1024C$.

\eqref{control-s} can be verified by using Lemma~\ref{lemma-beta}
and \eqref{s'''} and the fact that $|\Psi_{s''',\,n,\,p}(x)|\le
\beta_{n2}(x)+\beta_{n3}(x)$. This completes the proof of the lemma.
\eop

\begin{lemma}\label{divergence} Assume $n-q_{\max}\to\infty$ as
$n\to\infty$. Then as $n\to\infty$
\begin{equation}\label{div1}
\frac{(n-q_{\max})(n-p)}{p-q_{\max}}\log\big(1+\frac{p-q_{\max}}{3(n-p)}\big)>\frac13
\log(1+n-q_{\max})\to\infty
\end{equation}
and
\begin{equation}\label{div2}
(n-p)^2\log\big(1+\frac{p-q_{\max}}{3(n-p)}\big)\to\infty.
\end{equation}
In addition, if $p\to\infty$, then
\begin{equation}\label{div3}
p(n-p)\log\big(1+\frac{p-q_{\max}}{3(n-p)}\big)\to\infty.
\end{equation}
\end{lemma}

\noindent\textit{Proof.}  Define $f(x)=\frac{\log(1+x)}{x}$, $x>0$.
We can verify that
\[
f'(x)=-\frac{1}{x^2}\big(\log(1+x)+\frac{1}{1+x}-1\big)=-\frac1{x^2}\int^x_0\frac{t}{1+t^2}dt<0,
\]
that is, $f(x)$ is decreasing in $x>0$.   Now set
$x_n=(p-q_{\max})/3(n-p)$. Then $x_n< n-q_{\max}$. We get
\begin{eqnarray*}
\frac{(n-q_{\max})(n-p)}{p-q_{\max}}\log(1+\frac{p-q_{\max}}{3(n-p)})&=&\frac13(n-q_{\max})f(x_n)\\
&>& \frac13(n-q_{\max})f(n-q_{\max})\\
&=&\frac13\log(1+n-q_{\max})\\
&\to&\infty,
\end{eqnarray*}
proving \eqref{div1}.   \eqref{div2} follows from \eqref{div1} since
\[
\frac{(n-q_{\max})(n-p)}{p-q_{\max}}\Big/(n-p)^2=\frac{n-q_{\max}}{(p-q_{\max})(n-p)}=\frac1{p-q_{\max}}+\frac1{n-p}\le
2.
\]

One can easily verify that for any $y>0$, the function $x\log
(1+y/x)$ is increasing in $x\ge 1$. Therefore, we have
\[
p(n-p)\log\big(1+\frac{p-q_{\max}}{3(n-p)}\big)=\frac{p}{3}\times
3(n-p)\log\big(1+\frac{p-q_{\max}}{3(n-p)}\big)\ge
\frac{p}{3}\log\big(1+p-q_{\max}\big)\ge \frac{p}{3}\log(2)
\to\infty
\]
since $p\to\infty$ as $n\to\infty$. This proves \eqref{div3}. \eop

\begin{lemma}\label{Variance-Appr} If $n-q_{\max}\to\infty$ as
$n\to\infty$, then
\begin{equation}\label{small-small}
\frac{b(n,p)-\sum_{i=1}^kb(n, q_i)}{\beta_{n1}(0)}\to 0,
~~\frac{b(n,p)-\sum_{i=1}^kb(n, q_i)}{\sqrt{\beta_{n1}(0)}}\to 0,
\end{equation}
\begin{equation}\label{2-control}
\frac{p}{n(n-p)\sqrt{\beta_{n1}(0)}}\to
0,~~\frac{p}{n(n-p)^2\sqrt{\beta_{n1}(0)}}\to 0
\end{equation}
and
\begin{equation}\label{var-appr}
\lim_{n\to\infty}\frac{\mu_n-\bar\mu_n}{\sigma_n}=0,
~~~\lim_{n\to\infty}\frac{n^2\beta_{n1}(0)}{\sigma_n^2}=1
\end{equation}
as $n\to\infty$, where $b(n,q)$ is defined in \eqref{bnq}. Moreover,
\eqref{var-appr} implies \eqref{mean-var}.
\end{lemma}

\noindent\textit{Proof.}  We have $b(n,p)-\sum_{i=1}^kb(n,
q_i)=\frac14\beta_{n2}(0)$, where $\beta_{n2}$ is defined in
\eqref{beta-nr}. From \eqref{beta2-bound}, \eqref{beta1} and
\eqref{div1} we get
\[
\frac{b(n,p)-\sum_{i=1}^kb(n, q_i)}{\beta_{n1}(0)}=O\left(\Big(
\frac{(n-q_{\max})(n-p)}{p-q_{\max}}\log\big(1+\frac{p-q_{\max}}{3(n-p)}\big)\Big)^{-1}
\right)\to 0.
\]
Similarly,  from \eqref{beta2-bound}, \eqref{beta1} and \eqref{div2}
we have
\begin{eqnarray*}
\frac{b(n,p)-\sum_{i=1}^kb(n, q_i)}{\sqrt{\beta_{n1}(0)}}&=&O\Big(
\frac{\sqrt{p}(p-q_{\max})}{\sqrt{n}(n-q_{\max})(n-p)\sqrt{\log\big(1+\frac{p-q_{\max}}{3(n-p)}\big)}}\Big)\\
&=&O\Big(
\frac{1}{\sqrt{(n-p)^2\log\big(1+\frac{p-q_{\max}}{3(n-p)}\big)}}\Big)\\
&\to& 0.
\end{eqnarray*}
This completes the proof of \eqref{small-small}.

By using \eqref{beta1} and \eqref{div2} we get
\[
\frac{p}{n(n-p)\sqrt{\beta_{n1}(0)}}=O\big(\frac{\sqrt{p/n}}{\sqrt{(n-k)^2\log(1+\frac{p-q_{\max}}{3(n-p)})}}\big)
=O\big(\frac{1}{\sqrt{(n-k)^2\log(1+\frac{p-q_{\max}}{3(n-p)})}}\big)\to
0,
\]
and
\[
\frac{p}{n(n-p)^2\sqrt{\beta_{n1}(0)}}\le
\frac{p}{n(n-p)\sqrt{\beta_{n1}(0)}}\to 0,
\]
proving \eqref{2-control}.

The second limit in \eqref{var-appr} follows from
\eqref{small-small} since
$\sigma_n^2=n^2\beta_{n1}(0)+2n^2\big(b(n,p)-\sum_{i=1}^kb(n,
q_i)\big)$.

Set $\ell(x)=\log x-\frac{1}{2x}-\frac{1}{12x^2}$ for $x>0$.

It follows from \eqref{1storder} that
\[
g(x):=\psi(x)-\ell(x)=O(x^{-4})~~~~\mbox{ as }x\to\infty.
\]
Since $g(x)$ is a smooth function in $x>0$, it is easy to show that
\[
|g(x)|\le \frac{C}{x^3}, ~~~x\ge \frac14
\]
for some constant $C$.   Following the same lines in the proof of
\eqref{control-h} we have
\[
|\Psi_{g,\, n,\, p}(0)|\le \frac{Dp}{n(n-p)^2}
\]
for some $D>0$, which together with \eqref{2-control} implies
 \begin{equation}\label{diff}
\frac{n\Psi_{\psi,\, n,\, p}(0)-n\Psi_{\ell,\, n,\,
p}(0)}{\sqrt{n^2\beta_{n1}(0)}}=\frac{\Psi_{\psi,\, n,\,
p}(0)-\Psi_{\ell,\, n,\,
p}(0)}{\sqrt{\beta_{n1}(0)}}=\frac{\Psi_{g,\, n,\,
p}(0)}{\sqrt{\beta_{n1}(0)}}\to 0.
 \end{equation}

In virtue of \eqref{mean},
\begin{eqnarray*}
&&\sum^p_{j=1}\log(n-j)-\sum^k_{i=1}\sum^{q_i}_{j=1}\log(n-j)\\
&=&\sum^p_{j=1}\big(\log(n-j)-\log(n-1)\big)-\sum^k_{i=1}\sum^{q_i}_{j=1}
\big(\log(n-j)-\log(n-1)\big)\\
&=&(p-n+\frac12)\log(1-\frac{p}{n})-\sum^k_{i=1}(q_i-n+\frac12)\log(1-\frac{q_i}{n})+O(\frac{p}{n(n-p)}).
\end{eqnarray*}
We have used the condition $p=\sum^k_{i=1}q_i$ to simplify the
expression. We skip the details here. Similarly, by using
\eqref{var}, we obtain
\[
\sum_{j=1}^{p}\frac{1}{n-j}-\sum_{i=1}^{k}
\sum_{j=1}^{q_j}\frac{1}{n-j}=-\log(1-\frac{p}{n})+\sum^k_{i=1}\log(1-\frac{q_i}{n})+O(\frac{p}{n(n-p)}).
\]
Therefore, we have from \eqref{2-control} that
\begin{eqnarray*}
&&\Psi_{\ell,\, n,\, p}(0)\\
&=&\sum^p_{j=1}\log(n-j)-\sum^k_{i=1}\sum^{q_i}_{j=1}\log(n-j)-\big(\sum_{j=1}^{p}\frac{1}{n-j}-\sum_{i=1}^{k}
\sum_{j=1}^{q_j}\frac{1}{n-j}\big)-\frac{1}{3}\big(b(n,p)-\sum_{i=1}^kb(n,
q_i)\big)\\
&=&(p-n+\frac32)\log(1-\frac{p}{n})-\sum^k_{i=1}(q_i-n+\frac32)\log(1-\frac{q_i}{n})-\frac{1}{3}\big(b(n,p)-\sum_{i=1}^kb(n,
q_i)\big)+O(\frac{p}{n(n-p)})\\
&=&-\frac{\bar\mu_n}{n}+o\big(\sqrt{\beta_{n1}(0)}\big),
\end{eqnarray*}
that is,
\[
n\Psi_{\ell,\, n,\,
p}(0)=-\bar\mu_n+o\big(\sqrt{n^2\beta_{n1}(0)}\big),
\]
which coupled with \eqref{diff}, \eqref{var-appr} and the fact that
$\mu_n=-n\Psi_{\psi, \,n,\,p}(0)$ proves the first limit in
\eqref{var-appr}.

Finally, \eqref{mean-var} follows from \eqref{Variance-Appr} since
$\bar\sigma_n^2=n^2\beta_{n1}(0)$ in view of \eqref{bar-Var}. This
completes the proof of the lemma.  \eop

\section{Proofs of the main results}\label{proof}

\noindent\textit{Proof of Theorem~\ref{thm1}.}

We will prove \eqref{T0} under conditions $n-q_{\max}\to\infty$ and
$p\to\infty$.

Set $V_n=-2\log \Lambda_n$. Then it follows from \eqref{equ1} that
$V_n=-n\log W_n$ and
\[
\frac{-(V_n-\mu _n)}{\sigma_n}=\frac{n}{\sigma_n}\log
W_n+\frac{\mu_n}{\sigma_n}.
\]
In order to prove $\eqref{T0}$, it is sufficient to show the
moment-generating function of $-(V_n-\mu _n)/\sigma_n$ converges to
that of the standard normal, that is, for any $z$,
\begin{equation}\label{M-convergence}
\lim_{n\to\infty}\big(\log E(W_n^{nz/\sigma_n})+\frac{\mu_n
z}{\sigma_n}\big)= \frac{z^2}{2}.
\end{equation}

In view of \eqref{gamma}, \eqref{Delta} and \eqref{Psi}, we have for
$t>\frac{p-n}{2}$
\[
\log\frac{\Gamma_q(\frac{n-1}{2}+t)}{\Gamma_q(\frac{n-1}{2})}=\Delta_{
\log\Gamma,\,n,\,q}(t)-\Delta_{\log\Gamma, \,n, \,q}(0)
\]
and
\begin{eqnarray}\label{logM}
\log\Big(\mathbf{E}(W_n^t)\Big)&=&\log\Big(\frac{\Gamma_p(\frac{n-1}{2}+t)}{\Gamma_p(\frac{n-1}{2})}\prod_{i=1}^k
\frac{\Gamma_{q_i}(\frac{n-1}{2})}{\Gamma_{q_i}(\frac{n-1}{2}+t)}\Big)\nonumber\\
&=&\log\frac{\Gamma_p(\frac{n-1}{2}+t)}{\Gamma_p(\frac{n-1}{2})}-\sum^k_{i=1}
\log\frac{\Gamma_{q_i}(\frac{n-1}{2}+t)}{\Gamma_{q_i}(\frac{n-1}{2})}\nonumber\\
&=&\Psi_{\log\Gamma,\,n,\,p}(t)-\Psi_{\log\Gamma, \,n, \,p}(0).
\end{eqnarray}

Now we apply Taylor's theorem to expand $\Psi_{h,\,n,\, p}(t)$ with
a second order remainder. For any $t$ with $|t|\le (n-p)/4$, there
exists a $c_t$ with $|c_t|\le t\le (n-p)/4$ such that
\begin{equation}\label{h-expansion}
\Psi_{h, \,n, \,p}(t)-\Psi_{h, \,n, \,p}(0)=\Psi_{h', \,n,
\,p}(0)t+\Psi_{h'', \,n, \,p}(c_t)\frac{t^2}{2}.
\end{equation}
Note that we have used the property given in \eqref{Psi'}.
Similarly, we expand $\Psi_{h,\, n, \, p}(t)$ to the third order
\begin{equation}\label{s-expansion}
\Psi_{s, \,n, \,p}(t)-\Psi_{s, \,n, \,p}(0)=\Psi_{s', \,n,
\,p}(0)t+\Psi_{s'', \,n, \,p}(0)\frac{t^2}{2}+\Psi_{s''', \,n,
\,p}(\bar c_t)\frac{t^3}{6},
\end{equation}
where $|\bar c_t|\le |t|\le (n-p)/4$.

We will show
\begin{equation}\label{variance-bound}
    \lim_{n\to\infty}\frac{n^2}{\sigma_n^2(n-p)^2}=0.
\end{equation}
 Note that
\begin{equation}\label{var-beta}
\sigma_n^2\ge n^2\beta_{n1}(0).
\end{equation}
Then it follows from \eqref{beta1} that
\[
\sigma_n^2\ge \frac23np\log\big(1+\frac{p-q_{\max}}{3(n-p)}\big),
\]
which coupled with Lemma~\ref{divergence} implies that
\[
\frac{n^2}{\sigma^2_n(n-p)^2}\le
\frac{2n}{p(n-p)^2\log\big(1+\frac{p-q_{\max}}{3(n-p)}\big)}=\frac{2}{p(n-p)\log\big(1+\frac{p-q_{\max}}{3(n-p)}\big)}+
\frac{2}{(n-p)^2\log\big(1+\frac{p-q_{\max}}{3(n-p)}\big)}\to 0,
\]
proving \eqref{variance-bound}.

Now we proceed to prove \eqref{M-convergence} for any fixed $z$. For
fixed $z$,  set $t=t_n=\frac{nz}{\sigma_n}$.  Then it follows from
\eqref{variance-bound} that $t_n=o(n-p)$, which implies $|t_n|\le
(n-p)/4$ for all large $n$. Therefore, we can apply \eqref{logM},
\eqref{h-expansion} and \eqref{s-expansion} with $t=t_n$. From
Lemma~\ref{control} and \eqref{var-beta} we have
\[
|\Psi_{h'', \,n, \,p}(c_{t_n})|t_n^2\le
\frac{Dp}{n(n-p)^2}o((n-p)^2)=o(1)\]
 and
 \[
|\Psi_{s''', \,n, \,p}(c_{t_n})||t_n|^3\le
\frac{32\beta_{n1}(0)}{n-p}t_n^2o\big(n-p\big)=
o\big(\frac{n^2\beta_{n1}(0)}{\sigma_n^2}\big) =o(1).
\]
Then by combining \eqref{logM}, \eqref{hx}, \eqref{h-expansion} and
\eqref{s-expansion} that
\begin{eqnarray*}
\log\Big(\mathbf{E}(W_n^{nz/\sigma_n})\Big)
&=&\Psi_{\log\Gamma,\,n,\,p}(t_n)-\Psi_{\log\Gamma, \,n, \,p}(0)\\
&=&\Psi_{s, \,n, \,p}(t_n)+\Psi_{h, \,n, \,p}(t_n)-\big(\Psi_{s,
\,n, \,p}(0)+\Psi_{h,
\,n, \,p}(0)\big)\\
&=&\Psi_{s, \,n, \,p}(t_n)-\Psi_{s, \,n, \,p}(0)+\Psi_{h, \,n,
\,p}(t_n)-\Psi_{h,
\,n, \,p}(0)\\
&=&\big(\Psi_{s', \,n, \,p}(0)+\Psi_{h', \,n, \,p}(0)\big)t_n
+\Psi_{s'', \,n, \,p}(0)\frac{t_n^2}{2}+\Psi_{s''', \,n,
\,p}(c_{t_n})\frac{t_n^3}{6}+ \Psi_{h'', \,n,
\,p}(c_{t_n})\frac{t_n^2}{2}\\
 &=&\Psi_{\psi, \,n,
\,p}(0)t_n +\Psi_{s'', \,n, \,p}(0)\frac{t_n^2}{2}+o(1)\\
 &=&\frac{n\Psi_{\psi, \,n,
\,p}(0)}{\sigma_n}z +\frac{n^2\Psi_{s'', \,n, \,p}(0)}{\sigma_n^2}\frac{z^2}{2}+o(1)\\
 &=&-\frac{\mu_nz}{\sigma_n} +\frac{z^2}{2}+o(1),
\end{eqnarray*}
where we have used the following facts
\[
n\Psi_{\psi, \,n, \,p}(0)=-\mu_n, ~~~ n^2\Psi_{s'', \,n,
\,p}(0)=\sigma_n^2.
\]
This proves \eqref{M-convergence}.   \eop

\noindent\textit{Proof of Theorem~\ref{thm2}.}  Without loss of
generality, assume $q_1=q_{\max}$. When $n-q_{\max}=r$ and $n-p=v$
are fixed integers, $k$ is bounded and $q_i$, $2\le i\le k$ are also
bounded.  We can employ the subsequence argument to prove
\eqref{non-normal}, that is, for any subsequence of $n$, we will
show that there exists its further subsequence along which
\eqref{non-normal} holds.  Our criterion for selection of
subsequences is first to choose a subsequence of the given
subsequence of $n$, along which $k=k_n$ converges to a finite limit,
which implies $k_n$ is ultimately a constant, and then to select its
further subsequence along which $q_2$ converges. We repeat the same
procedure until we find a subsequence along which $q_k$ converges.
The last sub-sequence will be the one along which $k$ and $q_i$'s
are ultimately constant integers. The proof of \eqref{non-normal}
along such a subsequence is essentially the same as the proof when
$k=k_n$ is fixed and all $q_i$ for $2\le i\le k$ are also fixed
integers. For brevity, we will prove \eqref{non-normal} by assuming
$k$ and $q_i$ for $2\le i\le k$ are constants for all large $n$.

We first work on $W_n$ defined in \eqref{equ1}. From
Lemma~\ref{lemma1b}, $\log W_n-r\log n$ has the same distribution as
\[
\sum^k_{i=2}\sum^{q_i}_{j=1}\log V_{ij}-r\log
n=\sum_{i=2}^k\sum^{q_i}_{j=1}\log nV_{ij},
\]
where $q_i^*=\sum^{i-1}_{j=1}q_j$ for $2\le i\le k$, the $V_{ij}$'s
are independent random variables and $V_{ij}$ has a
beta($\frac12(n-q_i^*-j)$, $\frac12q_i^*$) distribution.

Set $\bar{q}_i=q_i^*-q_1=q_i^*-(p-r)=q_i^*-n+v+r$ for $2\le i\le k$.
We have $\bar{q}_2=0$, and $\bar{q}_i=\sum^{i-1}_{j=2}q_j$ if $i>2$.
This implies $\bar{q}_i$'s are fixed integers for all large $n$. For
any $1\le j\le q_i$, $2\le i\le k$,  $V_{ij}$ has a
beta\big($\frac12(r+v-\bar{q}_i-j)$, $\frac12(n-r-v+\bar{q}_i)$\big)
distribution.

It follows from Section 8.5 in Blitzstein and Hwang~\cite{BH2014}
that a beta($a$, $b$) random variable has the same distribution as
$X(a)/(X(a)+Y(b))$, where $X(a)$ and $Y(b)$ are two independent
random variables, $X(a)$ has a Gamma($a$, $2$) distribution with
density function $f(x)=x^{a-1}e^{-x/2}/2^a\Gamma(a)$, $x>0$, and
$Y(b)$ has a Gamma($b$, $2$) distribution.

Now for each pair of $(i,j)$ with $1\le j\le q_i$ and $2\le i\le k$,
set $a=\frac12(r+v-\bar{q}_i-j)$ and $b=\frac12(n-r-v+\bar{q}_i)$.
Then $a+b=\frac{1}{2}(n-j)$. Note that
$X(a)=X(\frac12(r+v-\bar{q}_i-j))$ and $Y(
b)=Y(\frac12(n-r-v+\bar{q}_i))$ are independent chi-square random
variables with $r+v-\bar{q}_i-j$ and $n-r-v+\bar{q}_i$ degrees of
freedom, respectively, and
$X(\frac12(r+v-\bar{q}_i-j))+Y(\frac12(n-r-v+\bar{q}_i))$ is also a
chi-square random variable with $n-j$ degrees of freedom. By using
the law of large numbers,
\[
\frac{X(\frac12(r+v-\bar{q}_i-j))+Y(\frac12(n-r-v+\bar{q}_i))}{n}=\frac{n-j}{n}\frac{X(\frac12(r+v-\bar{q}_i-j))+Y(\frac12(n-r-v+\bar{q}_i))}{n-j}
\]
converges in probability to $1$ as $n\to\infty$. Since $nV_{ij}$ has
the same distribution as
\[
X(\frac12(r+v-\bar{q}_i-j))\frac{n}{X(\frac12(r+v-\bar{q}_i-j))+Y(\frac12(n-r-v+\bar{q}_i))}
\]
which converges in distribution to a chi-square random variable with
$r+v-\bar{q}_i-j$ degrees of freedom, that is
\[
nV_{ij}\td Y_{r+v-\bar{q}_i-j},
\]
we have
\[
\sum^k_{i=2}\sum^{q_i}_{j=1}\log V_{ij}-r\log
n=\sum_{i=2}^k\sum^{q_i}_{j=1}\log nV_{ij}\td
\sum_{i=2}^k\sum^{q_i}_{j=1}\log
Y_{r+v-\bar{q}_i-j}=\sum^{r+v-1}_{j=v}\log Y_j,
\]
which is the limiting distribution of $\log W_n-r\log n$. We obtain
\eqref{non-normal} by noting that
\[
\frac{-2\log\Lambda_n+rn\log n}{n}=-\big(\log W_n-r\log n\big).
\]
This completes the proof of Theorem~\ref{thm2}. \eop

\noindent\textit{Proof of Theorem~\ref{thm3}.}

The sufficiency follows from Theorem~\ref{thm1}, that is, under
conditions $p\to\infty$ and $n-q_{\max}\to\infty$ as $n\to\infty$,
the central limit theorem \eqref{clt} holds with $a_n=\mu_n$ and
$b_n=\sigma_n$.

Now assume \eqref{clt} holds.  We need to show  $p\to\infty$ and
$n-q_{\max}\to\infty$.  If any one of the two conditions is not
true, there must exist a subsequence of $\{n\}$, say $\{n'\}$, along which\\
\noindent\textbf{a.} $p$ is fixed, $k$ is fixed and all $q_i$'s are
fixed, or\\
\noindent\textbf{b.} $p-q_{\max}=r$ and $n-p=v$ for some fixed
integers $r\ge 1$ and $v\ge1$.

Condition \textbf{b} holds when $n-q_{\max}$ is bounded because both
$q-q_{\max}$ and $n-p$ are bounded.

The subsequence $\{n'\}$ along which condition \textbf{a} holds can
be embedded in an entire sequence along which condition \textbf{a}
holds. Since the limiting distribution of $-2\log\Lambda_n$ is a
chi-square distribution according to \eqref{classic},  its
subsequential limit along $\{n'\}$ cannot be normal.     For the
same reason,  under condition \textbf{b},  the subsequential limit
is also non-normal from Theorem~\ref{thm2}; See Remark 4.  Under
either condition
 \textbf{a} or condition \textbf{b}, it results in a contradiction to the central limit theorem in
\eqref{clt}.  This completes the proof of the necessity.    \eop

\noindent\textit{Proof of Theorem~\ref{thm4}.}

Similar to the proof of Theorem 2 in Qi et al.~\cite{QWZ19},  we use
the subsequence argument to prove the theorem. It suffices to proved
\eqref{chisquare} under each of the following two assumptions:
\newline
\noindent{Case 1:} $p_n=p$ and $k_n=k$ and all $q_i$'s are fixed
integers for all large $n$;\\
\noindent{Case 2:} $p_n\to\infty$ and $n-q_{\max}\to\infty$ as
$n\to\infty$.

Under Case 1, $f_n=f$ is a constant for all large $n$, and
\eqref{classic} holds. Since $\rho_n$ defined in \eqref{rho}
converges to one,  $-2\log \Lambda_n$ converges in distribution to a
chi-square distribution with $f$ degrees of freedom.  Note that
$Z_n$ is defined in \eqref{zn}. To prove \eqref{chisquare},  it
suffices to verify that
\begin{equation}\label{case1-small}
\lim_{n\to\infty}\frac{2f}{\sigma_n^2}=1\mbox{ and
}\lim_{n\to\infty}\big(f-\mu_n\sqrt{\frac{2f}{\sigma_n^2}}\big)=0.
\end{equation}

Using the notation in the proof of Lemma~\ref{lemma-beta}, we have
$q_i^*$'s are fixed integers.  it follows form  \eqref{beta1=} that
\[
\beta_{n1}(0)=\frac{2(1+o(1))}{n^2}
\sum_{i=2}^{k}\sum^{q_i}_{j=1}q_i^*=\frac{2(1+o(1))}{n^2}
\sum_{i=2}^{k}q_iq_i^*=\frac{1+o(1)}{n^2}\big(p^2-\sum^k_{i=1}q_i^2\big)=\frac{2(1+o(1))f}{n^2},
\]
which together with \eqref{var-appr} implies $2f/\sigma_n^2\to 1$
and proves the first limit in \eqref{case1-small}. To prove the
second limit, it suffices to show $\lim_{n\to\infty}\mu_n=f$ or
equivalently $\lim_{n\to\infty}\bar\mu_n=f$ by using
\eqref{mean-var}.

From \eqref{var-appr} and \eqref{small-small},
$n\big(b(n,p)-\sum^k_{i=1}b(n,
q_i)\big)=o(n\beta_{n1}(0))=o(\frac1n)$. Then by using Taylor's
expansion we have from \eqref{mu1} that
\begin{eqnarray*}
\bar\mu_n&=&n\sum_{i=1}^k(q_i-n+\frac32)\log(1-\frac{q_i}{n})-n(p-n+\frac32)\log(1-\frac{p}{n})+o(\frac{1}{n})\\
&=&n\Big(\sum^k_{i=1}(q_i-n+\frac{3}{2})(-\frac{q_i}{n}-\frac12(\frac{q_i}{n})^2+O(\frac1{n^3}))- (p-n+\frac{3}{2})(\frac{p}{n}+\frac12(\frac{p}{n})^2+O(\frac1{n^3}))\Big)
+o(\frac1n)\\
&=&\frac12(p^2-\sum^k_{i=1}q_i^2)+O(\frac1n)\\
&=&f+o(1),
\end{eqnarray*}
This proves the second limit in \eqref{case1-small}.

The proof under Case 2 is the same as that in Qi et
al.~\cite{QWZ19}, and it is outlined as follows. First, rewrite
\eqref{chisquare} as
\begin{equation}\label{chisquare1}
\lim_{n\to\infty}\sup_{-\infty<x<\infty}|P(\frac{Z_n-f_n}{\sqrt{2f_n}}\le
x)-P(\frac{\chi^2_{f_n}-f_n}{\sqrt{2f_n}}\le x)|=0.
\end{equation}
We can show $f_n\to\infty$ under assumption
$\lim_{n\to\infty}(n-q_{\max})=\infty$. Since $\chi^2_{f_n}$ can be
written as a sum of $f_n$ independent chi-square random variables
with one degree of freedom, we have  from the central limit theorem
that
\[
\lim_{n\to\infty}\sup_{-\infty<x<\infty}|P(\frac{\chi^2_{f_n}-f_n}{\sqrt{2f_n}}\le
x)-\Phi(x)|=0.
\]
To show \eqref{chisquare1}, it suffices to prove that
\[
\frac{Z_n-f_n}{\sqrt{2f_n}} \td N(0,1),
\]
 which is a direct consequence of Theorem~\ref{thm1} since
$\frac{Z_n-f_n}{\sqrt{2f_n}}
=\frac{-2\log\Lambda_n-\mu_n}{\sigma_n}.$ This completes the proof
of Theorem~\ref{thm4}.  \eop

\null\hskip 20pt

\noindent\textbf{Acknowledgements.} The authors would like to thank
the two referees whose constructive suggestions have led to
improvement in the readability of the paper. The research of
Yongcheng Qi was supported in part by NSF Grant DMS-1916014.


\baselineskip 12pt
\def\ref{\par\noindent\hangindent 25pt}

\end{document}